\documentclass[12pt,leqno]{article}
\usepackage{amsfonts}
\usepackage{amsmath, amsthm, amsfonts, amssymb, color}
\usepackage{mathrsfs}
\usepackage{color}
\theoremstyle{definition}

\newcommand{\scr}[1]{\mathscr #1}
\definecolor{wco}{rgb}{0.5,0.2,0.3}

\numberwithin{equation}{section} \theoremstyle{remark}

\newcommand{\ua}{\uparrow}

\title{{\bf      Identifying    Constant Curvature Manifolds, Einstein Manifolds, and Ricci Parallel Manifolds }\footnote{Supported in
 part by  NNSFC (11771326, 11431014).} }
\author{
{\bf     Feng-Yu Wang  }\\
\footnotesize{ Center of Applied Mathematics, Tianjin University, Tianjin 300072, China}\\
 \footnotesize{ Department of Mathematics,
Swansea University, Singleton Park, SA2 8PP, United Kingdom}\\
\footnotesize{  wangfy@tju.edu.cn, F.-Y.Wang@swansea.ac.uk}}
\begin{document}
\allowdisplaybreaks
\def\R{\mathbb R}  \def\ff{\frac} \def\ss{\sqrt} \def\B{\mathbf
B}
\def\N{\mathbb N} \def\kk{\kappa} \def\m{{\bf m}}
\def\ee{\varepsilon}\def\ddd{D^*}
\def\dd{\delta} \def\DD{\Delta} \def\vv{\varepsilon} \def\rr{\rho}
\def\<{\langle} \def\>{\rangle} \def\GG{\Gamma} \def\gg{\gamma}
  \def\nn{\nabla} \def\pp{\partial} \def\E{\mathbb E}
\def\d{\text{\rm{d}}} \def\bb{\beta} \def\aa{\alpha} \def\D{\scr D}
  \def\si{\sigma} \def\ess{\text{\rm{ess}}}
\def\beg{\begin} \def\beq{\begin{equation}}  \def\F{\scr F}
\def\Ric{\mathcal Ric} \def\Hess{\text{\rm{Hess}}}
\def\e{\text{\rm{e}}} \def\ua{\underline a} \def\OO{\Omega}  \def\oo{\omega}
 \def\tt{\tilde}
\def\cut{\text{\rm{cut}}} \def\P{\mathbb P} \def\ifn{I_n(f^{\bigotimes n})}
\def\C{\scr C}      \def\aaa{\mathbf{r}}     \def\r{r}
\def\gap{\text{\rm{gap}}} \def\prr{\pi_{{\bf m},\varrho}}  \def\r{\mathbf r}
\def\Z{\mathbb Z} \def\vrr{\varrho} \def\ll{\lambda}
\def\L{\scr L}\def\Tt{\tt} \def\TT{\tt}\def\II{\mathbb I}
\def\i{{\rm in}}\def\Sect{{\rm Sect}}  \def\H{\mathbb H}
\def\M{\scr M}\def\Q{\mathbb Q} \def\texto{\text{o}} \def\LL{\Lambda}
\def\Rank{{\rm Rank}} \def\B{\scr B} \def\i{{\rm i}} \def\HR{\hat{\R}^d}
\def\to{\rightarrow}\def\l{\ell}\def\iint{\int}
\def\EE{\scr E}\def\Cut{{\rm Cut}}
\def\A{\scr A} \def\Lip{{\rm Lip}}\def\S{\mathbb S}
\def\BB{\scr B}\def\Ent{{\rm Ent}} \def\i{{\rm i}}\def\itparallel{{\it\parallel}}
\def\g{{\mathbf g}}\def\Sect{{\mathcal Sec}}\def\T{\mathcal T}
\maketitle

\begin{abstract} We establish  variational formulas for Ricci upper and lower bounds, as well as a  derivative formula for the Ricci curvature. Combining these with derivative and Hessian formulas of the heat semigroup  developed from stochastic analysis, we identify  constant curvature manifolds,  Einstein manifolds and Ricci parallel manifolds  by using  analytic  formulas and semigroup inequalities.
Moreover,    explicit Hessian estimates are derived for the heat semigroup on Einstein and Ricci parallel manifolds.
\end{abstract} \noindent
 \vskip 2cm

\section{Introduction}

Let $(M,\g)$ be a $d$-dimensional complete Riemannian manifold. Let $TM=\cup_{x\in M}T_xM$ be the bundle of tangent vectors of $M$, where for every $x\in M$, $T_xM$ is the tangent space at point $x$. We will denote $\<u,v\>=\g(u,v)$ and $|u|= \ss{\<u,u\>}$ for $u,v\in T_xM, x\in M.$ Let $\mathcal R$,   $\Ric$ and $\Sect$ denote  the Riemannian curvature tensor,   Ricci curvature and sectional curvature respectively.

When $M$ has constant Ricci curvature, i.e. $\Ric= K\g$ (simply denote by $\Ric=K$) for some constant $K$, the metric is a vacuum solution of Einstein field equations (in particular for $d=4$ in general relativity),  and $M$ is called an Einstein manifold.   In differential geometry, a basic problem is to characterize topological and geometry properties of Einstein manifolds. For instance, according to Thorpe \cite{Tho} and Hitchin \cite{Hi}, if a four-dimensional compact manifold admits an Einstein metric, then  $$\chi(M)\ge \ff 3 2 |\tau(M)|,$$ where $\chi$ is the Euler characterization and $\tau$ is the signature, and the equality holds if and only if $M$ is  a  flat torus, a Calabi-Yau manifold, or a quotient thereof. 
Recently, Brandle \cite{[4]}   showed that Einstein manifolds with positive isotropy curvature are space forms, while with nonnegative isotropy curvature are locally symmetric, and an ongoing  investigation is  to classify Einstein manifolds of both positive and negative sectional curvatures.

Since the curvature tensor is determined by the sectional curvature,  $M$ is called a constant curvature manifold if $\Sect=k$ for some constant $k\in\R$.  Complete simply connected constant curvature manifolds are called space forms, which are classified by  hyperbolic space (negative constant sectional curvature), Euclidean space (zero sectional curvature) and  unit sphere (positive sectional curvature) respectively.   All other connected complete constant curvature manifolds are quotients of space forms by some group of isometries.  Constant curvature manifolds are Einstein but not vice versa.

A slightly larger  class of  manifolds are Ricci parallel manifolds, where the Ricci curvature is constant under parallel transports; that is, $\nn \Ric=0$ where $\nn$ is the Levi-Civita connection. An Einstein manifold is Ricci parallel  but the inverse is not true. When the manifold is simply connected and indecomposable (i.e. does not split as non-trivial Riemannian products), these two properties are equivalent,   see e.g. \cite[Chapter XI]{KN}.

In this paper, we aim to identify the above three classes of manifolds by using integral formulas with respect to the volume measure and derivative inequalities of the heat semigroup. To state the main results, we introduce some notations where most are standard in the literature.

For $f,g\in C^2(M)$ and $x\in M$, consider the Hilbert-Schmidt inner product of the Hessian tensors $\Hess_f$ and $\Hess_g$:
$$\<\Hess_f,\Hess_g\>_{HS}(x)= \sum_{i,j=1}^d \Hess_f(\Phi^i,\Phi^j) \Hess_g(\Phi^i,\Phi^j),$$ where $\Phi=(\Phi^i)_{1\le i\le d}\in O_x(M),$ the space of orthonormal bases of   $T_xM$. Then the Hilbert-Schmidt norm of $\Hess_f$ reads
$$\|\Hess_f\|_{HS}= \ss{\<\Hess_f, \Hess_f\>_{HS}}\,.$$

 For a symmetric 2-tensor $T$ and a constant $K$, we write $\T\ge K$ if
$$\T(u,u)\ge K|u|^2,\ \ u\in TM.$$ Similarly,   $\T\le K$ means $\T(u,u)\le  K|u|^2, u\in TM.$   Let $\T^\#: TM\to TM$ be defined by
$$\<\T^\#(u),v\>= \T(u,v),\ \ u,v\in T_xM, x\in M.$$ Then $\T^\#$ is a symmetric map, i.e. $\<\T^\#u,v\>=\<\T^\#v,u\>$ for $u,v\in T_xM, x\in M.$ Let
$$\|\T\|(x)= \sup\big\{|\T^\#(u)|:\ u\in T_xM, |u|\le 1\big\},\ \ x\in M.$$ A $C^1$  map  $Q:TM\to TM$ with $QT_xM\subset T_x M$ for $x\in M$ is called  constant, if
$\nn (Q v)(x)=0$ holds for any $x\in M$ and vector field $v$ with  $\nn v(x)=0$. So,   $M$  is Ricci parallel if only if    $\Ric^\#$ is  a constant map.

For any symmetric $2$-tensor $\T$, define
\beq\label{RT} (\mathcal R\T)(v_1,v_2) := {\rm tr} \big\<\mathcal R(\cdot,v_2)v_1, \T^\#(\cdot)\big\> =\sum_{i=1}^d \big\<\mathcal R(\Phi^i, v_2)v_1, \T^\#(\Phi^i)\big\>,\end{equation}
where $v_1,v_2\in T_xM, x\in M, \Phi=(\Phi^i)_{1\le i\le d}\in O_x(M).$ Since $T$ is symmetric, so is $\mathcal RT$. Let
\beg{align*}&\|\mathcal R\|(x)= \sup \big\{\|\mathcal R \T\|(x):\ \T\ \text{is\ a\ symmetric\ 2-tensor},\ \|\T\|(x)\le 1\big\},\\
&\|\mathcal R\|_\infty = \sup_{x\in M} \|\mathcal R\|(x).\end{align*}

For a smooth tensor $\T$, consider the Bochner Laplacian
$$\DD \T:= {\rm tr} \big(\nn_\cdot\nn_\cdot \T\big).$$
Then $\ff 1 2\DD$ generates a contraction  semigroup $P_t =\e^{\ff t 2 \DD} $ in the $L^2$ space of tensors, see  \cite[Theorems 2.4 and 3.7]{Stri} for details.
In Subsection 3.1,  we will  prove a probabilistic formula of $P_t\T$, which is a smooth tensor when $\T$ is smooth with compact support.
Precisely, for any $x\in M$ and $\Phi\in O_x(M)$, let $\Phi_t(x)$ be the horizontal Brownian motion starting at $\Phi$, and let $X_t(x):=\pi \Phi_t(x)$ be the Brownian motion starting at $x$, see \eqref{E1} and \eqref{E2} below for details. Then $\parallel_t= \Phi_t(x)\Phi^{-1}: T_xM\to T_{X_t(x)}M$ is called the stochastic parallel transport along the Brownian path. Both $X_t(x)$ and $\parallel_t$ do not depend on the choice of the initial value $\Phi\in O_x(M).$ When the manifold is stochastically complete (i.e. the Brownian motion is non-explosive),  for a bounded $n$-tensor $\T$ we have
$$(P_t\T)(v_1,\cdots,v_n)= \E \big[\T (\parallel_t v_1,\cdots, \parallel_t v_n)\big],\ \ v_1,\cdots, v_n\in T_xM.$$
 In the following, we will take this regular (rather than $L^2$) version of the heat semigroup $P_t$.

Finally,  For $v\in T_xM$, let  $W_t(v)\in T_{X_t(x)}M$ solve  the following covariant differential equation
\beq\label{WT}\ff{\d }{\d t} \Phi_t(x)^{-1}W_t(v) = - \ff 1 2\Phi_t^{-1}(x)\Ric^\#(W_t(v)),\ \ W_0(v)=v.\end{equation}
$W_t$ is called the damped stochastic parallel transport. When  $\Ric\ge K$ for some constant $K\in \R$, we have $|W_t(v)|\le \e^{-\ff K 2 t}|v|,\ t\ge 0, v\in T_xM.$

In Section 2 and Sections 6-8, we will present a number of identifications of   constant curvature manifolds, Einstein manifolds, and Ricci parallel manifolds.
In particular, the following assertions are direct consequences of    Theorems \ref{T1.0}, \ref{TN}, \ref{T1.1} and \ref{T1.3} below.

\paragraph{(A)  Constant curvature.} Let $k\in \R$. Each of the following assertions is equivalent to $\Sect =k$: \beg{enumerate}
\item[$(A_1)$]  For any $ t\ge 0$ and $f\in C_0^\infty(M),$
$ \Hess_{P_t f}= \e^{-d kt} P_t\Hess_f + \ff 1 d (1-\e^{-dkt}) (P_t\DD f)\g.$
\item[$(A_2)$]  For any   $f\in C_0^\infty(M),$
$\Hess_{\DD f}-   \DD \Hess_f= 2k \big\{(\DD f)\g-d\,  \Hess_f\big\}.$
\item[$(A_3)$]  For any   $f\in C_0^\infty(M),$
$$\ff 12 \DD \|\Hess_f\|_{HS}^2 -\<\Hess_{\DD f}, \Hess_f\>_{HS} -\|\nn \Hess_f\|_{HS}^2 = 2k \big(d\|\Hess_f\|_{HS}^2 -(\DD f)^2\big),$$
where $\|\nn \Hess_f\|_{HS}^2:= \sum_{i=1}^d \|\nn_{\Phi^i} \Hess_f\|_{HS}^2,\ \Phi=(\Phi^i)_{1\le i\le d} \in O_x(M).$
\item[$(A_4)$]  For any  $x\in M, u\in T_xM$ and  $f\in C_0^\infty(M)$ with $\Hess_f(x)= u\otimes u\ ($i.e. $\Hess_f(v_1,v_2)= \<u,v_1\>\<u,v_2\>, v_1,v_2\in T_xM)$,
    $$\big(\Hess_{\DD f} - \DD\Hess_f\big)(v,v)= 2k \big(|u|^2|v|^2-\<u,v\>^2\big),\ \ v\in T_xM.$$
   \end{enumerate}

According to the Bochner-Weitzenb\"ock formula, for any constant $K\in\R$,   $\Ric=K$ is equivalent to each of the following formulas:
\beg{align*} &\nn P_tf= \e^{-\ff t 2 K}P_t\nn f,\ \ f\in C_0^\infty(M),\ \ t\ge 0,\\
&\DD\nn f-\nn \DD f =K\nn f,\ \ f\in C_0^\infty(M),\\
&\ff 1 2 \DD |\nn f|^2  -\<\nn\DD f, \nn f\> - \|\Hess_f\|_{HS}^2 = K|\nn f|^2,\ \ f\in C_0^\infty(M). \end{align*}
So, $(A_1)$-$(A_3)$ can be regarded as the corresponding formulas for $\Sect=k.$ Below, we present some other identifications of Einstein manifolds.

\paragraph{(B)  Einstein manifolds.} Let $\mu$ be the volume measure, and denote $\mu(f)=\int_Mf\d\mu$ for $f\in L^1(\mu)$. $M$ is Einstein if and only if
\beq\label{EIN1}\beg{split} &\mu\big(\<\nn f,\nn g\>\big) \mu\big((\DD f)^2-\|\Hess_f\|_{HS}^2\big) \\
&= \mu\big( |\nn f|^2\big)\mu\big((\DD f)(\DD g)-\<\Hess_f,\Hess_g\>_{HS}\big),\ \ f,g\in C_0^\infty(M).\end{split} \end{equation}
Moreover, for any $K\in \R$, $\Ric=K$ is equivalent to each of  the following assertions:
\beg{enumerate} \item[$(B_1)$]   $ \|\mathcal R\|_\infty <\infty$, and for any $ t\ge 0, f\in C_0^\infty(M),$
$$\Hess_{P_tf}= \e^{-Kt}P_t\Hess_f +\int_0^t \e^{-Ks} P_s(\mathcal R \Hess_{P_{t-s}f})\d s.$$
\item[$(B_2)$] For any $f\in C_0^\infty(M)$,
 $$\ff 1 2 \big\{\Hess_{\DD f}- \DD\Hess_f\big\}= (\mathcal R\Hess_f) - K\Hess_f.$$
\item[$(B_3)$] For any $f\in C_0^\infty(M)$,
$$ \mu\big((\DD f)^2-\|\Hess_f\|_{HS}^2\big) =K\mu\big(|\nn f|^2\big).$$
\item[$(B_4)$] For any $f,g\in C_0^\infty(M)$,
$$ \mu\big((\DD f)(\DD g)-\<\Hess_f,\Hess_g\>_{HS}\big) =K\mu\big(\<\nn f,\nn g\>\big).$$
\item[$(B_5)$]  There exists   $h:  [0,\infty)\times M\to   [0,\infty)$ with $\lim_{t\to 0} h(t,\cdot)=0$ such that
   $$ \big|P_t|\nn f|^2 - \e^{Kt} |\nn P_t f|^2 \big|\le h(t,\cdot) \big(\|\Hess_{P_tf}\|_{HS}^2+P_t\|\Hess_f\|_{HS}^2\big),\ \ t\ge 0, f\in C_0^\infty(M).$$
\end{enumerate}

\paragraph{(C)  Ricci Parallel manifolds.}    $M$ is Ricci parallel if and only if
$$\int_M \big\{(\DD f)^2- \|\Hess_f\|_{HS}^2 -\<Q\nn f,\nn f\>\big\}\d\mu =0,\ \ f\in C_0^\infty(M) $$ holds for some constant symmetric linear map $Q: TM\to TM$, and in this case $\Ric^\#=Q$. They are also equivalent to
each of the following statements:
\beg{enumerate} \item[$(C_1)$]   There exists a function
 $h: [0,\infty)\times M \to  [0,\infty) $ with $\lim_{t\downarrow 0}t^{-\ff 1 2}h(t,\cdot)=0$ such that
$$\|\Hess_{P_tf} -   P_t\Hess_f\|\le h(t,\cdot)\big( P_t\|\Hess_f\|+\|\Hess_{P_tf}\|\big),\ \ t\ge 0, f\in C_b^2(M).$$
\item[$(C_2)$]  $ \|\mathcal R\|_\infty <\infty$, and for any $x\in M, t\ge 0$, $f\in C_0(M)$ and $v_1,v_2\in T_xM$,
\beg{align*}    \Hess_{P_tf}(v_1,v_2) -  \E\big[\Hess_f(W_t(v_1), W_t( v_2))\big]
  = \E\int_0^t  \big(\mathcal R \Hess_{P_{t-s}f}\big)\big(W_{s}(v_1), W_{s}(v_2)\big) \d s.\end{align*}
  \item[$(C_3)$] For any $f\in C_0^\infty(M)$ and $x\in M$,
  $$\big(\Hess_{\DD f}- \DD \Hess_f\big)(v_1,v_2)= 2\big(\mathcal R\Hess_f\big)(v_1,v_2) - 2\Ric(v_1,\Hess_f^\#(v_2)),\ \ v_1,v_2\in T_xM.$$
  \item[$(C_4)$] For any $x\in M$ and  $f\in C_0^\infty(M)$  with $\Hess_f(x)=0$,
$$(\DD\Hess_f)(v_1,v_2)=\Hess_{\DD f}(v_1,v_2),\ \   v_1,v_2\in T_xM.$$
 \end{enumerate}
Since a symmetric $2$-tensor is determined by its diagonal,
one may take $v_1=v_2$ in $(C_2)$-$(C_4)$.

\

To prove these results, we establish the following variational and derivative formulas for $\Ric$, see Remark 2.1 and Theorem  \ref{T1.2} below.

\paragraph{(D) Formulas of $\Ric$.} Let $\overline{\Ric}$ and $\underline{\Ric}$ be the exact upper and exact lower bounds of $\Ric$. Then
\beg{align*} &\underline{\Ric}=\inf\Big\{\mu\big( (\DD f)^2-\|\Hess_f\|_{HS}^2\big):\ f\in C_0^\infty(M), \mu(|\nn f |^2)=1\Big\},\\
&\overline{\Ric}=\sup\Big\{\mu\big( (\DD f)^2-\|\Hess_f\|_{HS}^2\big):\ f\in C_0^\infty(M), \mu(|\nn f |^2)=1\Big\}.\end{align*}
Moreover, let $x\in M$, $v_1,v_2\in T_xM$. For any $f\in C_b^4(M)$ with $\nn f(x)=v_1, \Hess_f(x)=0$,
$$ (\nn_{v_2}\Ric)(v_1,v_1) =2 \lim_{t\downarrow 0} \ff{(P_t\Hess_f-\Hess_{P_tf})(v_1,v_2) }{t}= \big(\DD\Hess_f -\Hess_{\DD f}\big)(v_1,v_2). $$

The remainder of the paper is organized as follows. In Section 2, we  present variational formulas of  Bakry-Emery-Ricci upper and lower bounds, as well as integral characterizations of Einstein and Ricci parallel manifolds. In Section 3, we recall
   derivative  and Hessian formulas of  $P_t$ developed from stochastic analysis. Using these formulas we estimate the Hessian of $P_t$ in Section 4, and establish
  a formula for $\nn\Ric$ in Section 5. Finally,
by applying results presented in Sections 3-5, we identify constant curvature manifolds, Einstein manifolds, and Ricci parallel manifolds in Sections  6-8 respectively.

 \section{Characterizations of Bakry-Emery-Ricci curvature}

Let $V\in C^2(M)$,   $\mu_V(\d x)= \e^{V(x)}\mu(\d x)$, and $L_V=\DD+\nn V$. Then $L_V$ is symmetric in $L^2(\mu_V).$
Consider the Bakry-Emery-Ricci curvature
$$\Ric_V:= \Ric -\Hess_V.$$
\beg{defn}   The manifold $M$ is called $V$-Einstein if $\Ric_V=K$ for some constant $K\in\R$, while it is called $\Ric_V$ parallel if $\nn \Ric_V=0 \ ($i.e. $\Ric_V^\#: TM\to TM$ is a constant map). \end{defn}

By the integral formula of Bochner-Weitzenb\"ock, we have
\beq\label{BWI} \int_M\big\{(L_Vf)^2 -\|\Hess_f\|_{HS}^2-\Ric_V(\nn f,\nn f)\big\}\d\mu_V=0,\ \ f\in C_0^\infty(M).\end{equation} According to Theorem \ref{T1.0} below,   this formula   identifies the curvature $\Ric_V$, and provides sharp upper and lower bounds of $\Ric_V$, as well as integral characterizations of $V$-Einstein and $\Ric_V$ parallel manifolds.

\beg{thm}\label{T1.0} Let $K\in\R$ be a constant, and let $Q: TM\to TM$ be a symmetric continuous linear map. \beg{enumerate}
\item[$(1)$] $\Ric^\#_V=Q$ if and only if
\beq\label{LK0} \int_M \big\{(L_V f)^2 -\|\Hess_f\|_{HS}^2 -\<Q\nn f,\nn f\>\big\}   \d \mu_V =0,\ \ f\in C_0^\infty(M).\end{equation}Consequently,
 $M$ is $\Ric_V$ parallel if and only if $\eqref{LK0}$ holds for  some   symmetric  constant linear  map $Q:TM\to TM$.
\item[$(2)$] For any $V\in C^2(M)$ and $K\in \R$,
$\Ric_V\ge K$ if and only if
\beq\label{RL}  \int_M \Big\{(L_V f)^2-\|\Hess_f\|_{HS}^2\Big\}\d\mu_V\ge K \int_M|\nn f |^2\d\mu_V,\ \  f\in C_0^\infty(M);\end{equation}
while $\Ric_V\le K$ if and only if
\beq\label{RU}  \int_M \Big\{(L_V f)^2-\|\Hess_f\|_{HS}^2\Big\}\d\mu_V\le K \int_M|\nn f |^2\d\mu_V,\ \  f\in C_0^\infty(M).\end{equation}
\item[$(3)$]  $M$ is   $V$-Einstein   if and only if
\beq\label{00}\beg{split} &\mu_V\big(|\nn f|^2\big)\cdot\mu_V\big((L_V f)(L_V g) - \<\Hess_f,\Hess_g\>_{HS} \big)\\
&=\mu_V\big( \<\nn f,\nn g\>\big)\cdot\mu_V\big( (L_V f)^2-\|\Hess_f\|_{HS}^2\big),\ \  f,g\in C_0^\infty(M).\end{split}\end{equation}
Moreover, for any constant $K\in \R$,   $\Ric_V=K$ is equivalent to each of
\beq\label{01}\int_M \big\{(L_V f)^2 -\|\Hess_f\|_{HS}^2 -K|\nn f|^2\big\}  \d \mu_V =0,\ \ f\in C_0^\infty(M),\end{equation}
and \beq\label{01'}\int_M \big\{(L_V f)(L_V g) -\<\Hess_f,\Hess_g\>_{HS}   -K\<\nn f,\nn g\> \big\}  \d \mu_V =0,\ \ f,g\in C_0^\infty(M).\end{equation}
\end{enumerate}
\end{thm}

\paragraph{Remark 2.1.} Theorem \ref{T1.0}(2) provides the following  variational formulas of the upper and lower bounds of $\Ric_V$.
Let
$$\overline{\Ric_V}=\sup\{\Ric_V(u,u): \ u\in TM, |u|=1\},\ \ \underline{\Ric_V}=\inf\{\Ric_V(u,u): \ u\in TM, |u|=1\}.$$ We have
\beg{align*} &\underline{\Ric_V}=\inf\Big\{\mu_V\big( (L_V f)^2-\|\Hess_f\|_{HS}^2\big):\ f\in C_0^\infty(M), \mu_V(|\nn f |^2)=1\Big\},\\
&\overline{\Ric_V}=\sup\Big\{\mu_V\big( (L_V f)^2-\|\Hess_f\|_{HS}^2\big):\ f\in C_0^\infty(M), \mu_V(|\nn f |^2)=1\Big\}.\end{align*}

To prove   Theorem \ref{T1.0}, we need the following lemma.

\beg{lem}\label{LP2} For a continuous symmetric linear map $Q: TM\to TM$, if
\beq\label{NG} \int_M \<Q\nn f,\nn f\>\d\mu_V\ge 0, \ \ f\in C_0^\infty(M),\end{equation}
then $Q\ge 0$; that is, $\<Qu,u\>\ge 0$ for $u\in TM.$\end{lem}

\beg{proof} Using $\e^VQ$ replacing $Q$, we may and do assume that $V=0$ so that $\mu_V=\mu$ is the volume measure.

 {\bf (a)} We first consider $M=\R^d$ for which we have $Q=(q_{ij})_{1\le i,j\le d}$ with $q_{ij}=q_{ji}$ for  some continuous functions  $q_{ij}$ on $\R^d$. Thus,
$$\<Qu,v\>= \sum_{i,j=1}^d q_{ij} u_iv_j,\ \ u,v\in \R^d.$$
Without loss of generality, we only prove that $Q(0)\ge 0$. Using the eigenbasis of $Q(0)$, we may and do assume that $Q(0)={\rm diag}\{q_1,\cdots, q_d\}.$ It suffices to prove $q_l\ge 0$ for $1\le l\le d.$ For any $f\in C_0^\infty(\R^d)$, let
$$f_n(x)= f\circ \phi_n(x),\ \ (\phi_n(x))_i:= n x_i\ {\rm if}\ i\ne l, \ (\phi_n(x))_l:= n^2 x_l,\ \ n\ge 1. $$
By \eqref{NG},
\beg{align*} 0 &\le \sum_{i,j=1}^d \int_{\R^d} q_{ij}(x) (\pp_i f_n)(x) \pp_j f_n(x) \d x \\
&= \int_{\R^d} \Big(n^4 q_{ll}(x) (\pp_l f)^2\circ\phi_n(x) + 2n^3 \sum_{j\ne l}  q_{jl}(x) \{(\pp_lf)(\pp_j f)\}\circ\phi_n(x) \\
&\qquad\qquad  +n^2\sum_{i,j\ne l} q_{ij}(x) \{(\pp_if)(\pp_j f)\}\circ\phi_n(x) \Big)\d x\\
&= \int_{\R^d} \Big(n^{3-d} q_{ll}\circ \phi_n^{-1}(x) (\pp_l f)^2 (x) + 2n^{2-d} \sum_{j\ne l}  q_{jl}\circ \phi_n^{-1}(x) \{(\pp_lf)(\pp_j f)\} (x)\\
&\qquad \qquad +n^{1-d}\sum_{i,j\ne l} q_{ij}\circ \phi_n^{-1}(x) \{(\pp_if)(\pp_j f)\} (x) \Big)\d x,\end{align*} where the last step is due to the integral transform $x\mapsto \phi_n^{-1}(x)$.
Since $\phi_n^{-1}(x)\to 0$ as $n\to\infty$, multiplying both sides by $n^{d-3}$ and letting $n\to\infty$ we arrive at
$$\int_{\R^d}  q_l (\pp_lf)^2(x) \d x \ge 0,\ \ f\in C_0^\infty(\R^d).$$ Thus, $q_l\ge 0$ as wanted.

{\bf (b)} In general, for $x_0\in M$, we take a neighborhood $O(x_0)$ of $x_0$ such that it is  diffeomorphic to $\R^d$ with $x_0$ corresponding to $0\in\R^d$. Let $\psi: O(x_0)\to \R^d$ with $\psi(x_0)=0$ be a diffeomorphism. Then under the local charts induced by $\psi$,
$$\<Q\nn(f\circ\psi),\nn(f\circ \psi)\>\d\mu= \sum_{i,j=1}^d q_{ij}(x) (\pp_if)(x)(\pp_jf) (x)\d x,\ \ f\in C_0^\infty(\R^d) $$ holds for some symmetric matrix-valued  continuous functional $(q_{ij})_{1\le i,j\le d}$.
Therefore, by step {\bf (a)}, \eqref{NG} implies $\sum_{i,j=1}^d q_{ij}(x) u_iu_j\ge 0, x,u\in \R^d.$ In particular, $Q(x_0)\ge 0.$  \end{proof}

 \beg{proof} [Proof of Theorem \ref{T1.0}]  (1)   By the Bochner-Weitzenb\"ock formula, we have
\beq\label{BWQ}   \ff 1 2 L_V |\nn f|^2 -\<\nn f,\nn L_V f\>= \Ric_V(\nn f,\nn f) +\|\Hess_f\|_{HS}^2,\ \ f\in C_0^\infty(M).\end{equation}
So, $\Ric^\#_V=Q$ implies
$$\ff 1 2 L_V |\nn f|^2 -\<\nn f,\nn L_V f\>= \<Q\nn f,\nn f\> +\|\Hess_f\|_{HS}^2.$$ Integrating both sides with respect to $\mu_V$  proves
\eqref{LK0}.

On the other hand,
 integrating both of \eqref{BWQ} with respect to $\mu_V$, we obtain \eqref{BWI}.   This together with \eqref{LK0} implies
\beq\label{LK} \int_{\R^d} \big\<\Ric^\#_V(\nn f)-Q\nn f,\nn f\big\>\,\d \mu_V=0,\ \ f\in C_0^\infty(M). \end{equation} Therefore, by Lemma \ref{LP2} for $\Ric_V^\#-Q$ replacing $Q$, we prove $\Ric_V^\#=Q$.

 (2) By integrating both sides of \eqref{BWQ} with respect to $\mu_V$, we see that  $\Ric_V^\#\ge K$ implies \eqref{RL}.
  On the other hand, applying  Lemma \ref{LP2} to $Q=\Ric_V^\#-K$, if \eqref{RL} holds then  $\Ric^\#_V\ge K.$ Similarly, we can    prove the equivalence of $\Ric^\#_V\le K$  and \eqref{RU}. Therefore, Theorem \ref{T1.0}(2) holds.

  (3) By assertion (2),   $\Ric_V=K$ is equivalent \eqref{01}.  It remains to prove that \eqref{00} is equivalent to the Einstein property, since this together with the equivalence of $\Ric_V=K$ and \eqref{01} implies the equivalence of \eqref{01} and \eqref{01'}.

 If $M$ is $V$-Einstein, there exists a constant $K\in\R$ such that $\Ric_V= K$. Let $f,g\in C_0^\infty(M)$ with $\mu_V(|\nn f|^2)>0$, let $f_s=f+sg$. Then there exists $s_0>0$ such that $\mu_V(|\nn f_s|^2)>0$ for  $s\in [0,s_0]$. By \eqref{01}, we have
 $$h(s):= \ff{\mu_V((L_V f_s)^2 -\|\Hess_{f_s}\|_{HS}^2)}{\mu_V(|\nn f_s|^2)}= K,\ \ s\in [0,s_0].$$ So,
 \beg{align*}&\ff{\mu_V(|\nn f|^2) \mu_V((L_V f)(\DD g)- \<\Hess_f,\Hess_g\>) - \mu_V(\<\nn f,\nn g\>) \mu_V((L_V f)^2-\|\Hess_f\|_{HS}^2)}{\mu_V(|\nn f|^2)^2} \\
 & = \ff 1 2 h'(0)=0.\end{align*}
 Therefore, \eqref{00} holds.

 On the other hand, for $f\in C_0^\infty(M)$ with $\mu_V(|\nn f|^2)>0$, let
 $$K= \ff{\mu_V((L_V f)^2-\|\Hess_f\|_{HS}^2)}{\mu_V(|\nn f|^2)}\in \R.$$
 By the equivalence of $\Ric^\#_V=K$ and \eqref{01},   it suffices to prove
  \beq\label{ER} \mu_V\big((L_V g)^2-\|\Hess_g\|_{HS}^2\big) = K \mu_V(|\nn g|^2),\ \ g\in C_0^\infty(M).\end{equation}
  By the definition of $K$, this formula holds when $g$ is a linear combination of $f$ and $1$. So, we assume that $f, g$ and $1$ are linear independent. In this case,
  $$g_s:= (1-s)f+sg,\ \ s\in [0,1]$$ satisfies $\mu_V(|\nn g_s|^2)>0.$ Let
  $$K(s):= \ff{\mu_V((L_V g_s)^2 -\|\Hess_{g_s}\|_{HS}^2)}{\mu_V(|\nn g_s|^2)},\ \ s\in [0,1].$$ Then   \eqref{00} implies
  \beg{align*} K'(s)= \ff 2 {\mu_V(|\nn g_s|^2)^2}& \Big\{\mu_V(|\nn g_s|^2) \mu_V\big((L_Vg_s)L_V(g-f)- \<\Hess_{g_s}, \Hess_{g-f}\>_{HS}\big)\\
   &- \mu_V(\<\nn g_s, \nn(g-f)\>) \mu_V\big((L_V g_s)^2-\|\Hess_{g_s}\|_{HS}^2\big)\Big\}=0,\ \ s\in [0,1].\end{align*}
  Therefore,
  $$\ff{\mu_V((L_V g)^2-\|\Hess_g\|_{HS}^2)}{\mu_V(|\nn g|^2)} = K(1) = K(0)=K,$$ that is, \eqref{ER} holds as desired.
 \end{proof}

\section{Derivative and Hessian formulas of $P_t$ }

In this section, by using the (horizontal) Brownian motion, we first formulate the heat semigroup $P_t$ acting on tensors,  then  recall the derivative and Hessian  formulas of   $P_t$   on functions.

\subsection{ Brownian motion and heat semigroup on tensors}

Consider the projection operator  from   the orthonormal frame bundle $O(M)$ onto $M$:
$$\pi: O(M)\to M;\ \ \pi\Phi =x\ \text{if}\ \Phi\in O_x(M).$$ Then for any $a\in \R^d$ and $\Phi=(\Phi^i)_{1\le i\le d}\in O(M)$,
$$\Phi a :=\sum_{i=1}^d a_i\Phi^i\in T_{\pi\Phi }M;$$
and  for any $v\in T_{\pi\Phi}M$,
$$\Phi^{-1} v := \sum_{i=1}^d \<v,\Phi^i\> e_i\in \R^d,$$
where $\{e_i\}_{1\le i\le d}$ is the canonical orthonormal basis of $\R^d$.

For any $a\in\R^d$ and $\Phi\in O(M)$, let $\Phi(s)$ be the parallel transport of $\Phi$ along the geodesic $s\mapsto \exp[s\Phi a], s\ge 0.$ We have $$H_a(\Phi):=\ff{\d}{\d s} \Phi(s)\big|_{s=0} \in T_\Phi O(M),$$ which is a horizontal vector field on $O(M)$. Let $H_i=H_{e_i}$. Then
$(H_i)_{1\le i\le d}$ forms the canonical  orthonormal basis for the space of horizontal vector fields:
$$H_a= \sum_{i=1}^d a_i H_i,\ \ a=(a_i)_{1\le i\le d}\in \R^d.$$
We call
\beq\label{E11}\DD_{O(M)}:= \sum_{i=1}^d H_i^2\end{equation}  the horizontal Laplacian on $O(M).$ For any smooth $n$-tensor $\T$,
$$\T^{O(M)}(\Phi):= \big(\T(\Phi^{i_1},\cdots, \Phi^{i_n})\big)_{1\le i_1,\cdots, i_n\le d}\in \otimes^n\R^d,\ \ \Phi=(\Phi^i)_{1\le i\le d}\in O(M)$$ gives rise to a smooth $ \otimes^n\R^d$-valued function on $O(M)$. Let  $\DD$ be the Bochner Laplacian $\DD$.
For any $x\in M$, $v_1,\cdots, v_n\in T_xM$ and $\Phi\in O_x(M),$ we have
\beq\label{*ZJ}\beg{split} &(\DD \T)(v_1,\cdots, v_n) =  \{\DD_{O(M)} \T^{O(M)}(\Phi)\}(\Phi^{-1}v_1,\cdot,\Phi^{-1}v_n),\\
&(\nn_v\T)(v_1,\cdots, v_n) = (\nn_{H_{\Phi^{-1}v}} \T^{O(M)})(\Phi^{-1} v_1,\cdots, \Phi^{-1} v_n),\ \ v\in T_xM.\end{split} \end{equation}

Now, consider the following SDE on $O(M)$:
\beq\label{E1} \d\Phi_t= H(\Phi_t)\circ\d B_t :=\sum_{i=1}^d H_i(\Phi_t)\circ \d B_t^i,\end{equation} where $B_t:=(B_t^i)_{1\le i\le d}$ is the $d$-dimensional Brownian motion on a complete probability space $(\OO, \F, \P)$ with natural filtration $\F_t^B:=\si(B_s:s\le t)$ (by convention, we always take the completion of a $\si$ field).   The solution is called the Horizonal Brownian motion on $O(M)$. Let $\zeta$ be the life time of the solution. When $\zeta=\infty$ (i.e. the solution is non-explosive) we call the manifold $M$ stochastically complete. It is the case  when
\beq\label{RC} \Ric\ge - c (1+\rr^2) \end{equation} for some constant $c>0$, where $\rr$ is the Riemannian distance to some fixed point, see the proof of Proposition \ref{PP} blow. See also   \cite{Hsu} and references within for the stochastic completeness under weaker conditions.

Let $X_t:=\pi\Phi_t$ for $t\in [0,\zeta)$. Then $(X_t)_{t\in [0,\zeta)} $    solves
the SDE
\beq\label{E2} \d X_t = \Phi_t(X_t)\circ \d B_t,\end{equation} and is called the Brownian motion on $M$.   For any  $x\in M$, let $X_t(x)$ be the solution of \eqref{E2} with $X_0=x$. Note that   $X_t(x)$ does not depend on the choice of the initial value $\Phi_0\in O_x(M)$ for \eqref{E1}, and for fixed initial value  $\Phi_0\in O_x(M)$,
both $\Phi_t(X_t)$ and $B_t$ are measurable with  respect to $\F_t^X:=\si(X_s: s\le t).$ Therefore, $\F_t^B=\F_t^X, t\ge 0.$
When the   initial point $x$ is clearly given in the context,   we will simply denote $X_t(x)$ by  $X_t$.

Let
$${\it \parallel}_t:= \Phi_t\Phi_0^{-1}: T_{X_0}M\to T_{X_t}M $$  be the stochastic parallel transport along $X_t$, which also does not depend on the choice of $\Phi_0$.
For any  smooth $n$-tensor $\T$ with compact support, let
\beq\label{ZJJ} (P_t\T)(v_1,\cdots,v_n)= \E \big[1_{\{t<\zeta\}}\T (\parallel_t v_1,\cdots, \parallel_t v_n)\big],\ \ v_1,\cdots, v_n\in T_xM.\end{equation}  As is well known in the function (i.e. $0$-tensor) setting, by \eqref{E11}, the first formula in \eqref{*ZJ}, \eqref{E1}  and It\^o's formula
we have the forward/backward Kolmogorov  equations
\beq\label{DS5}  \pp_t P_t \T = \ff 1 2 \DD P_t \T = \ff 1 2 P_t \DD\T, \ \
 \T\in C_0^\infty. \end{equation}

For later use, we present the following exponential estimate of the Brownian motion under condition \eqref{RC}.
\beg{prp}\label{PP} Assume $\eqref{RC}$. Then there exist  constants $r, c(r)>0$ such that
 \beq\label{SB}\E \exp[\e^{-r t}\rr^2(X_t)]\le \exp\Big[\rr^2(x)+\ff{c(r)(1-\e^{-rt})}r \Big],\ \ t\ge 0, x\in M.\end{equation} Consequently, if
\beq\label{AS4}\lim_{\rr\to \infty} \ff{\log \|\nn\Ric\|}{\rr^2}=0,\end{equation}
  then  for any $\vv\in (0,1)$ and $p\ge 1$, there exists a positive function $C_{\vv,p}\in C([0,\infty))$ such that
\beq\label{AS3} \E \|\nn \Ric\|^p(X_t(x)) \le \e^{\vv\rr^2(x)+ C_{\vv,p} (t)},\ \ t\ge 0, x\in M.\end{equation}
\end{prp}
\beg{proof}  Let $\rr$ be the Riemannian distance to a fixed point $o\in M$. By    the Laplacian comparison theorem, \eqref{RC}  implies
$$\DD\rr\le c_1(\rr+\rr^{-1})$$   outside $\{o\}\cup {\rm cut}(o)$ for some constant $c_1>0$, where ${\rm cut}(o)$ is the cut-locus of $o$. By  It\^o's formula of $\rr(X_t)$ given in \cite{Ken}, this gives
$$\d\rr(X_t)^2\le c_2\{1 + \rr (X_t)^2\}\d t +2\rr(X_t)\d b_t$$ for some constant $c_2>0$ and an one-dimensional Brownian motion $b_t$. By It\^o's formula, for any $r>2+c_2$ there exists a constant $c(r)>0$ such that
\beg{align*} \d \exp[\e^{-r t}\rr^2(X_t) ] &\le \exp[\e^{-r t}\rr^2(X_t)]\Big(\big\{(c_2+c_2\rr^2(X_t))\e^{-rt} - r \e^{-rt}\rr^2(X_t)\\
 &\qquad \qquad \qquad \qquad + 2 \e^{-2rt} \rr^2(X_t)\big\}\d t + 2\e^{-rt} \rr(X_t) \d b_t\Big)\\
&\le \exp[\e^{-r t}\rr^2(X_t)]\big\{c(r)\e^{-rt} \d t + 2\e^{-rt} \rr(X_t) \d b_t\big\}.\end{align*}
Therefore, \eqref{SB} holds.

 On the other hand, by  \eqref{AS4} we may find a positive function $c_{p/\vv}\in C([0,\infty))$ such that
$$ \|\nn\Ric^\#\|^{p/\vv}\le \exp\big[\e^{-t} \rr^2 + c_{p/\vv}(t)\big],\ \ t\ge 0.$$ Combining this with \eqref{SB} for $r=1$, we obtain
\beg{align*} &\E \|\nn \Ric\|^p(X_t(x)) \le \big(\E \|\nn \Ric\|^{p/\vv}(X_t)\big)^\vv \\
&\le \big(\e^{\rr^2(x)+ c(1) +c_{\vv/p}(t) }\big)^\vv \le \e^{\vv \rr^2(x)+\vv c(1)+\vv c_{p/\vv}(t)}.\end{align*} Therefore, \eqref{AS3} holds for $C_{\vv,p}(t):= \vv [c(1)+c_{p/\vv}(t)].$\end{proof}

\subsection{Derivative formula of $P_t$}

 In this subsection, we assume
 \beq\label{CC2}  \Ric\ge -h(\rr) \ \text{ for\ some\ positive} \ h\in C([0,\infty))\ \text{with}\  \lim_{r\to\infty} \ff{h(r)}{r^2}=0.\end{equation}
Let $W_t: T_xM\to T_{X_t}M$ be defined in \eqref{WT}. By \eqref{CC2} and Proposition \ref{PP}, we have
\beq\label{SB3}\E \sup_{s\in [0,t]} \|W_s\|^p <\infty,\ \ t>0.\end{equation}
  We have (see e.g. \cite{EL, Hsu})
\beq\label{E4}\nn_v P_t f(x)=  \E \< \nn f(X_t(x)), W_t(v)\>,\ \ t\ge 0, f\in C_b^1(M).\end{equation}

Following the idea of \cite{EL, Th}, it is standard to establish the Bismut type formula using \eqref{E4}. By $\pp_t P_t f= \ff 1 2\DD P_t f$, \eqref{E2} and It\^o's formula, we have
$$\d P_{t-s}f(X_s)   = \<\nn P_{t-s} f(X_s), \Phi_s\d B_s\>,$$ so that
\beq\label{QB}f(X_t(x))= P_t f(x) +\int_0^t \<\nn P_{t-s} f(X_s(x)), \Phi_s(x)\d B_s\>.\end{equation}
In particular,  $P_{t-s}f(X_s)$ is a martingale. Next, by  \eqref{E4} and the Markov property,  we have
 $$\big\<\nn P_{t-s}f(X_s), W_s(v)\big\>=\E\big(\big\<\nn f(X_t), W_t(v)\big\> \big|\F_s^B\big),\ \ s\in [0,t],$$ which is again a martingale. Indeed, according to  e.g. \cite[(2.2.7)]{W14}, we have
\beq\label{BQ} \d\<\nn P_{t-s} f (X_s), W_s(v)\>= \Hess_{P_{t-s}f}(\Phi_s\d B_s, W_s(v)),\ \ s\in [0,t].\end{equation}
So, for any adapted process $h\in C^1([0,t])$ such that $h_0=0, h_t=1$,   \eqref{E4}, \eqref{QB} and \eqref{BQ} imply
\beg{align*} &\E\bigg[f(X_t(x)) \int_0^t \dot h_s \<W_s(v), \Phi_s(x)\d B_s\>\bigg] = \E\int_0^t \dot h_s \<W_s(v), \nn P_{t-s} f(X_s(x))\>\d s\\
&=\E\bigg[h_s\<W_s(v), \nn P_{t-s}f(X_s(x))\>\big|_0^t -\int_0^t h_s\d\<W_s(v),\nn P_{t-s}f(X_s(x))\>\bigg]\\
&= \E \<W_t(v), \nn f(X_t(x))\>= \nn_v P_t f(x).\end{align*} Therefore,
\beq\label{E10}\nn_v P_t f(x) =\E\bigg[f(X_t(x))\int_0^t \dot h_s \<W_s(v), \Phi_s(x)\d B_s\>\bigg],\ \ t>0, x\in M, f\in C_b^1(M).\end{equation}
This type formula is named after J.-M. Bismut, K.D.  Elworthy and X.-M. Li  because of  their  pioneering work \cite{Bis} and  \cite{EL}. The present version is due to \cite{Th} and has been  applied in \cite{AT, TW}   to derive gradient estimates  using local curvature conditions.

\subsection{Hessian formula of $P_t$}

To calculate $\Hess_{P_tf}$, we   introduce the following doubled  damped parallel transport $W_t(v_1,v_2)$ for $v_1, v_2\in T_xM$:
\beq\label{W2} \beg{split} \Phi_t(x)^{-1} W_t^{(2)}(v_1,v_2)= &\ff 1 2 \int_0^t\Phi_s^{-1}\big\{(\tt\nn \Ric^\#)(W_s(v_2))W_s(v_1)-\Ric^\#(W_s^{(2)}(v_1,v_2))\big\}\d s \\
&\qquad +   \int_0^t  \Phi_s(x)^{-1} \mathcal R (\Phi_s(x) \d B_s, W_s(v_2))W_s(v_1),\ \  t\ge 0,\end{split} \end{equation}
where the cycle derivative $\tt\nn \Ric^\#$ is defined by  \beq\label{W2'}\<(\tt\nn \Ric^\#)(u_2)u_1,u_3\>:= (\nn_{v_3} \Ric)(u_1,u_2)- (\nn_{u_1}\Ric) (u_2,u_3)- (\nn_{u_2}\Ric) (u_1,u_3) \end{equation} for $u_1,u_2,u_3\in T_yM, y\in M.$
According to Proposition \ref{PP} and \eqref{SB3},   conditions  \eqref{AS4} and \eqref{CC2} imply
$$\E \sup_{s\in [0,t]} \|W_s^{(2)}\|^p <\infty,\ \ p\ge 1, t>0.$$

\beg{prp}[\cite{APT03,Li}]\label{P1} Assume  $\eqref{AS4}$ and $\eqref{CC2}$. Then for any $f\in C_b^2(M)$ and $v_1,v_2\in T_xM$,
\beq\label{HESS}   \Hess_{P_tf}(v_1,v_2) = \E\big\{\Hess_f(W_t(v_1), W_t(v_2)) + \<\nn f(X_t(x)), W_t^{(2)}(v_1,v_2)\>\big\}.\end{equation}
\end{prp}
\beg{proof} Let $v_2(s)$ be the parallel transport of $v_2$ along the geodesic $s\mapsto \exp[s v_1], s\ge 0$. According to \eqref{AS3}, we define the following covariant derivative of $W_t$:
$$W^{(2)}_t(v_1,v_2) := \nn_{v_1} W_t(v_2) = \ff{\d}{\d s} W_t(v_s(s))\big|_{s=0}.$$
By \eqref{E4},
$$\Hess_{P_tf}(v_1,v_2)= \E\big\{\Hess_f(W_t(v_1), W_t(v_2)) + \<\nn f(X_t(x)), W^{(2)}_t(v_1,v_2)\>\big\}.$$ It remains to prove that $W_t^{(2)}$ satisfies \eqref{W2}. Since in the present setting  $\F_t^X=\F^B_t$, this follows from    formula (7) in \cite{APT03}, see also  (3.1) in \cite{Li}.
\end{proof}

In the same spirit of deducing the Bismut type derivative formula \eqref{E10}  from \eqref{E4},   Bismut type Hessian formulas of $P_tf$ have been presented in \cite{APT03, EL, Li} by using  \eqref{HESS}.
In Section 4 we will use the following local version of Hessian formula, which follows from \cite[Theorem 2.1]{APT03} and \cite[Proof of Theorem 3.1]{APT03} for e.g. $D_1= B(x,1), D_2= B(x, 2)$,
where $B(x,r)$ is the open geodesic ball at $x$ with radius $r$.

\beg{prp}[\cite{APT03}]\label{P2} Let $M$ be a complete noncompact Riemannian manifold.  Let
$$\tau_i(x)= \inf\{t\ge 0: X_t(x)\in \pp B(x, i)\},\ \ i=1,2.$$ There exists a positive function $C\in C(M)$ such that for any $x\in M$,
  $v_1,v_2\in T_xM$ with $|v_1,|v_2|\le 1$, and $f\in \B_b(M)$,
\beq\label{HU2} \Hess_{P_t f}(v_1,v_2)= \E \big[P_{t-t\land \tau_1(x)}f(X_{t\land\tau_1(x)}(x)) M_t+ P_{t-t\land\tau_2(x)}f(X_{t\land \tau_2(x)}(x))N_t\big],\ \ t> 0. \end{equation}
holds for some   adapted continuous processes $(M_t, N_t)_{t\ge 0}$ determined by $(X_t(x))_{0\le t\le \tau_2(x)}$ such that
\beq\label{HU1} \E\big[|N_{t}|+ |M_{t}|\big]\le \ff {C(x)}{t\land 1},\ \ t>0. \end{equation}

\end{prp}

\section{Hessian estimates   and  applications   }

In this section, we first present Hessian estimates of $P_t$ for Einstein and Ricci parallel manifolds, then apply these  results to describe the   lower and upper bounds  of the Ricci curvature.

Recall that for any $x\in M$ and $f\in C^2(M)$,
\beg{align*}&\|\Hess_f\|(x):= \sup\{|\Hess_f(u,v)|:\ u,v\in T_xM, |u|, |v|\le 1\},\\
&\|\Hess_f\|_{HS}^2(x):= \sum_{i,j=1}^d \Hess_f(\Phi^i,\Phi^j)^2,\ \ \Phi=(\Phi^i)_{1\le i\le d} \in O_x(M).\end{align*}

\beg{thm}\label{T3.1} Let $M$ be a Ricci parallel manifold with $\|\mathcal R\|_\infty<\infty$. Then for any $x\in M, t\ge 0$, $f\in C_0(M)$ and $v_1,v_2\in T_xM$,
\beq\label{HS5''}\beg{split}  &\Hess_{P_tf}(v_1,v_2) -  \E\big[\Hess_f(W_t(v_1), W_t( v_2))\big]\\
 &= \E\int_0^t \big(\mathcal R \Hess_{P_{t-s}f}\big)\big(W_{s}(v_1), W_{s}(v_2)\big) \d s,\end{split}\end{equation}where $\mathcal R\Hess_{P_{t-s}f}$ is defined in $\eqref{RT}$ for $T=\Hess_{P_{t-s}f}$. Consequently:
\beg{enumerate} \item[$(1)$] If $\Ric\ge K$, then for any $f\in C_b^2(M)$,
\beq\label{HS1} \|\Hess_{P_tf}\|\le \e^{( \|\mathcal R\|_\infty -K)t} P_t\|\Hess_f\|,\ \ t\ge 0.\end{equation}
 \item[$(2)$] If   $\Ric=K$, then
\beq\label{HS2} \|\Hess_{P_tf}\|_{HS}^2\le \e^{2( \|\mathcal R\|_\infty -K)t} P_t\|\Hess_f\|_{HS}^2,\ \ t\ge 0.\end{equation}\end{enumerate}
\end{thm}

\beg{proof} We fix $t>0$ and $ f\in C_b^2(M)$. Let ${\bf d} $ be the exterior differential. By e.g. \cite[(2.2.6)]{W14} we have
\beq\label{4.4*}\d ({\bf d} P_{t-s}f)(X_s)= \nn_{\Phi_s\d B_s} ({\bf d} P_{t-s}f)(X_s) +\ff 1 2 \Ric (\cdot, \nn P_{t-s}f(X_s))\d s, \ \ s\in [0,t].\end{equation} Equivalently,
\beq\label{HS3}\Phi_t^{-1} \nn f(X_t)= \Phi_0^{-1} \nn P_tf + \ff 1 2 \int_0^t \Phi_s^{-1}\Ric^\#(\nn P_{t-s}f(X_s))\d s+ \int_0^t\Hess_{P_{t-s}f}^\#(\Phi_s\d B_s).\end{equation}
On the other hand, since $\tt\nn \Ric^\#=0$, \eqref{W2} becomes
$$\Phi_t^{-1} W_t^{(2)}(v_1,v_2) = \int_0^t \Phi_s^{-1} \mathcal R (\Phi_s\d B_s, W_s(v_2)) W_s(v_1) -\ff 1 2 \int_0^t\Phi_s^{-1}\Ric^\#(W_s^{(2)}(v_1,v_2))\d s.$$
Combining this with \eqref{HS3*}, we obtain
\beq\label{HS4} \E\<\nn f(X_t), W_t^{(2)}(v_1, v_2)\> = \E\int_0^t {\rm tr} \Big\{\Hess_{P_{t-s}f}(\cdot, \mathcal R(\cdot, W_s(v_2))W_s(v_s))\Big\}. \end{equation} Plugging   \eqref{HS4} into \eqref{HESS} gives
 \beg{align*} &\Hess_{P_tf}(v_1,v_2)- \E\big[\Hess_f(W_t(v_1), W_t( v_2))\big] \\
 &=  \E\int_0^t{\rm tr} \big(
\big\< \mathcal R(\cdot, W_{s}(v_2))W_{s}(v_1), \Hess_{P_{t-s}f}^\#(\cdot)\big\>\big)\d s\\
&= \E\int_0^t  \big(
  \mathcal R\Hess_{P_{t-s}f}\big)(W_s(v_1), W_{s}(v_2)) \d s.\end{align*}   Therefore,    \eqref{HS5''} holds.

Below we prove \eqref{HS1} and \eqref{HS2} for Ricci parallel and Einstein manifolds respectively.

{\bf (a)} \eqref{HS5''} and $\Ric\ge K$ imply \eqref{HS1}. If $\Ric\ge K$, then \eqref{WT} implies
$$|W_t(v)|\le\e^{-\ff 1 2 Kt}|v|.$$ So,  according to \eqref{HS5''}, for any $s>0$ we have
$$\|\Hess_{P_sf}\|\le \e^{-Ks} P_s\|\Hess_f\| +  \|\mathcal R\|_\infty\int_0^s\e^{-Kr} P_r \|\Hess_{P_{s-r}f}\|\d r.$$
Letting
$$\phi(s)= \e^{-K(t-s)} P_{t-s} \|\Hess_{P_sf}\|,\ \ s\in [0,t],$$ we obtain
\beg{align*} \phi(s)&\le \e^{-K(t-s)} P_{t-s} \bigg(\e^{-sK} P_s\|\Hess_f\| +  \|\mathcal R\|_\infty \int_0^s \e^{-rK} P_r \|\Hess_{P_{s-r}f}\|\d r\bigg)\\
&\le \e^{-Kt} P_t\|\Hess_f\| +  \|\mathcal R\|_\infty \int_0^t \e^{-K(t+r-s)} P_{t+r-s} \|\Hess_{P_{s-r}f}\|\d r.\end{align*} Using the change of variable   $\theta=s-r$, we arrive at
$$\phi(s)\le \phi(0) + \|\mathcal R\|_\infty \int_0^s \phi(\theta)\d\theta,\ \ s\in [0,t].$$
By Gronwall's lemma, this implies
$$\phi(t)\le \phi(0)\e^{ \|\mathcal R\|_\infty t},$$ which is equivalent to \eqref{HS1}.

{\bf (b)} Let $\Ric=K$. Then \eqref{WT} implies $W_s(v)=\e^{-\ff K 2 s}v$. So, for $x\in M, v\in T_xM$ and $\Phi_0\in O_x(M)$,    \eqref{HS5''} implies
\beq\label{HS5'} \beg{split}\Hess_{P_tf}(v,\Phi_0^k) =\,& \e^{-Kt} \E\big[\Hess_f(\itparallel_{t} v, \Phi_t^k)\big]\\
& +    \E\int_0^t \e^{-sK}\big(\mathcal R
\Hess_{P_{t-s}f}\big) (\Phi_s^k,  \itparallel_{s}v)\d s.\end{split}\end{equation}
Let
$$\phi_{k} (s)= \e^{-K(t-s)}\, \E \|\Hess_{P_sf}^\#(\Phi_{t-s}^k)\|,\ \ s\in [0,t], 1\le k\le d.$$
By \eqref{HS5'} and the Markov property, for any $0\le s_2<s_1\le t$   we have
\beg{align*}&\Hess_{P_{s_1}f}(\itparallel_{t-s_1}v, \Phi_{t-s_1}^k) =    \e^{-K(s_1-s_2) }\, \E\big(\Hess_{P_{s_2}f}(\itparallel_{ t-s_2}v, \Phi_{t-s_2}^k)\big|\F_{t-s_1}^x\big)\\
&+   \int_0^{s_1-s_2}  \e^{-rK}\E\Big( \big(\mathcal R\Hess_{P_{s_1-r}f}\big) (\Phi_{t-s_1+r}^k,  \itparallel_{t-s_1+r}v)\Big|\F_{t-s_1}^x\Big)\d r.\end{align*}
So,
\beg{align*} &I_{k,v}(s_1,s_2)\\
&:= \E\Big|\e^{-(t-s_1)K}\Hess_{P_{s_1}f}(\itparallel_{t-s_1}v, \Phi_{t-s_1}^k)-   \e^{-K(t-s_2) }
 \E\big(\Hess_{P_{s_2}f}(\itparallel_{t-s_2}v, \Phi_{t-s_2}^k)\big|\F_{t-s_1}^x\big)\Big|\\
&\le  |v|\e^{-(t-s_1)K} \|\mathcal R\|_\infty   \E\int_0^{s_1-s_2}\e^{-rK} |\Hess_{P_{s_1-r}f}^\#(\Phi_{t-s_1+r}^k)|\d r\\
&= |v|\|\mathcal R\|_\infty  \int_{s_2}^{s_1} \e^{-K(t-\theta)}\E |\Hess_{P_\theta f}^\# (\Phi_{t-\theta}^k)| \d\theta\\
&= |v|\|\mathcal R\|_\infty \int_{s_1}^{s_2}  \phi_{k}(\theta)\d\theta,\end{align*}
where we have used the change of variable $\theta=s_1-r$.
Then
\beg{align*}&\phi_{k}(s_1)-\phi_{k}(s_2)  \le \sup_{|v|\le 1} I_{k,v}(s_1,s_2)\\
&\le  \|\mathcal R\|_\infty   \int_{s_1}^{s_2} \phi_{k}(\theta)\d\theta,\ \ 0\le s_2\le s_1\le t.\end{align*}
By Gronwall's lemma, this implies
$$|\Hess_{P_tf}^\#(\Phi_0^k)|=\phi_k(t)\le \e^{\|\mathcal R\|_\infty t}\phi_k(0)= \e^{(\|\mathcal R\|_\infty-K)t} \E  |\Hess_f^\#(\Phi_t^k)|,\ \ 1\le k\le d.$$ Therefore,
\beg{align*}   \|\Hess_{P_tf}\|^2_{HS}&=\sum_{k=1}^d |\Hess_{P_tf}^\#(\Phi_0^k)|^2 \le \e^{2(\|\mathcal R\|_\infty-K)t}\sum_{k=1}^d \big(\E  |\Hess_f^\#(\Phi_t^k)|\big)^2\\
&\le  \e^{2(\|\mathcal R\|_\infty-K)t}P_t   \|\Hess_f\|_{HS}^2. \end{align*}
\end{proof}

Next, we apply the above results to characterize the   lower and upper bounds of $\Ric$ for Ricci parallel manifolds.

\beg{thm}\label{TP3} Let $M$ be a Ricci parallel manifold. Then for any constant $K\in \R$, the following statements are equivalent each other:
\beg{enumerate}\item[$(1)$]  $\Ric\ge K.$
 \item[$(2)$]  For any $f\in C_0^\infty(M)$ and $t\ge 0$,
$$\ff{\e^{Kt}- \e^{2(K- \|\mathcal R\|_\infty)t}}{ 2 \|\mathcal R\|_\infty-K}   \|\Hess_{P_tf}\|^2  \le P_t|\nn f|^2 - \e^{Kt} |\nn P_t f|^2.$$
   \item[$(3)$] For any $f\in C_0^\infty(M)$ and $t\ge 0$,
$$\|\Hess_{P_tf}\|^2\int_0^t \ff{\e^{Ks}-\e^{2(K- \|\mathcal R\|_\infty)s}}{2 \|\mathcal R\|_\infty-K}\d s
 \le P_tf^2-(P_tf)^2-\ff{\e^{Kt}-1}K  |\nn P_t f|^2.$$
 \item[$(4)$]  For any $f\in C_0^\infty(M)$ and $t\ge 0$,
\beg{align*} & P_tf^2-(P_tf)^2- \ff{1-\e^{-Kt}}K  P_t |\nn f|^2 \\
& \le - \|\Hess_{P_tf}\|^2\,\e^{2(K- \|\mathcal R\|_\infty)t} \int_0^t \ff{\e^{2( \|\mathcal R\|_\infty-K)s}-\e^{-Ks}}{2 \|\mathcal R\|_\infty-K}\d s.\end{align*}\end{enumerate}
 \end{thm}
\beg{proof}{\bf  (a)} $(1)\Rightarrow (2)$. Let $t>0$ and $f\in C_b^2(M)$. By \eqref{E2} and It\^o's formula, we have
\beq\label{QBC} \beg{split}\d |\nn P_{t-s}f|^2(X_s) = \,&\bigg(\ff 1 2 \DD |\nn P_{t-s}f|^2(X_s)-\<\nn P_{t-s}f, \nn \DD P_{t-s}f\>(X_s)\bigg)\d s\\
&+2 \<\nn |\nn P_{t-s}f|^2(X_s), \Phi_s\d B_s\>,\ \ s\in [0,t].\end{split}\end{equation}
By the Bochner-Weitzenb\"ock formula and  $\Ric\ge K$,  we obtain
\beg{align*} &\ff 1 2 \DD |\nn P_{t-s}f|^2(X_s)-\<\nn P_{t-s}f, \nn \DD P_{t-s}f\>(X_s)\\&= \Ric(\nn P_{t-s}f, \nn P_{t-s}f)(X_s) +\|\Hess_{P_{t-s}f}\|_{HS}^2(X_s)\\
&\ge K |\nn P_{t-s}f|^2(X_s) + \|\Hess_{P_{t-s}f}\|_{HS}^2(X_s).\end{align*} Then \eqref{QBC} implies
$$\d |\nn P_{t-s}f|^2(X_s)\ge  \big(K |\nn P_{t-s}f|^2 + \|\Hess_{P_{t-s}f}\|_{HS}^2 \big)(X_s)\d s+
2 \<\nn |\nn P_{t-s}f|^2(X_s), \Phi_s\d B_s\>$$   for $s\in [0,t].$
Combining this with \eqref{HS1}, we arrive at
\beg{align*} &P_t |\nn f|^2 - \e^{Kt} |\nn P_t f|^2 \ge \int_0^t \e^{K(t-s)} P_s\|\Hess_{P_{t-s} f}\|_{HS}^2\d s\\
&\ge \int_0^t \e^{K(t-s)} \e^{-2( \|\mathcal R\|_\infty -K)s} \|\Hess_{P_t f}\|^2\d s\\
&= \ff{\e^{Kt}- \e^{2(K- \|\mathcal R\|_\infty)t}}{ 2 \|\mathcal R\|_\infty-K}    \|\Hess_{P_tf}\|^2 .\end{align*}

{\bf (b)} (2) implies  (3) and $(4).$ By \eqref{E2} and It\^o's formula, we have
$$\d (P_{t-s}f)^2(X_s) = |\nn P_{t-s} f|^2(X_s)\d s + \<\nn |P_{t-s}f|^2(X_s), \Phi_s\d B_s\>,\ \ s\in [0,t].$$
So,
\beq\label{PTT} P_t f^2- (P_tf)^2= \int_0^t P_s |\nn P_{t-s}f|^2 \d s.\end{equation}
Combining this with (2) and \eqref{HS1}, we obtain
\beg{align*} &P_t f^2-(P_tf)^2\\
&\ge \int_0^t \Big\{|\nn P_t f|^2\e^{Ks} + \ff{\e^{Ks}- \e^{2(K- \|\mathcal R\|_\infty)s}}{2 \|\mathcal R\|_\infty -K} \|\Hess_{P_tf} \|^2 \Big\}\d s.\end{align*} Then (3) is proved.

Similarly, \eqref{HS1} and (2) imply
\beg{align*}  &\e^{-Ks} P_t|\nn f|^2- P_{t-s}|\nn P_s f|^2\\
&\ge \ff{1-\e^{(K-2  \|\mathcal R\|_\infty)s}}{2 \|\mathcal R\|_\infty-K} P_{t-s}\|\Hess_{P_sf}\|^2\\
&\ge \e^{2(K- \|\mathcal R\|_\infty)t} \|\Hess_{P_tf}\|^2 \cdot \ff{\e^{2( \|\mathcal R\|_\infty-K)s}-\e^{-Ks}}{2 \|\mathcal R\|_\infty-K},
\end{align*}  which together with \eqref{PTT} gives (4).

{\bf (c)} Each of $(3)$ and $(4)$ implies $(1)$.  For  $v\in T_xM$ with $|v|=1$, take $f\in C_0^\infty(M)$ such that
$$\nn f(x)=v,\ \ \Hess_f(x)=0.$$ We have
\beq\label{LOP} \lim_{t\downarrow 0} \|\Hess_{P_tf}\|^2\int_0^t \ff{\e^{Ks}-\e^{2(K- \|\mathcal R\|_\infty)s}}{2 \|\mathcal R\|_\infty-K}\d s=0.\end{equation}
On the other hand,  by the Bochner-Weitzenb\"ock formula we have (see  \cite[Theorem 2.2.4]{W14}),
\beq\label{BW}\ff 1 2 \Ric(v,v)= \lim_{t\downarrow 0} \ff{P_t|\nn f|^p(x)- |\nn P_tf|^p(x)}{pt},\ \ p>0.\end{equation}
Combining this with $(3)$, \eqref{PTT} and \eqref{LOP}, we obtain
 \beg{align*} &0\le 2\lim_{t\downarrow 0} \ff {P_t f^2-(P_tf)^2 -\ff{\e^{Kt}-1}K |\nn P_tf|^2} {t^2} \\
&= \lim_{t\downarrow 0} \ff {2}{t^2} \int_0^t \big\{P_s|\nn P_{t-s}f|^2- \e^{Ks} |\nn P_t f|^2\big\}(x)\d s= \Ric(v,v)-K,\end{align*}
Therefore, $(3)$ implies $(1)$. Similarly, $(4)$ also implies $(1)$.
\end{proof}

The following result provides corresponding characterizations for the Ricci upper bound.

\beg{thm}\label{TP4} Let $M$ be a Ricci parallel manifold. Then for any constant $K\in \R$, the following are equivalent each other:
\beg{enumerate}\item[$(1)$]  $\Ric\le K.$
 \item[$(2)$]  For any $f\in C_0^\infty(M)$ and $t\ge 0$,
$$\ff{\e^{(2  \|\mathcal R\|_\infty-K) t}-1}{ 2 \|\mathcal R\|_\infty-K}  d P_t \|\Hess_{f}\|^2  \ge P_t|\nn f|^2 - \e^{Kt} |\nn P_t f|^2.$$
   \item[$(3)$] For any $f\in C_0^\infty(M)$ and $t\ge 0$,
\beg{align*}&dP_t \|\Hess_{f}\|^2\int_0^t \ff{\e^{(2 \|\mathcal R\|_\infty-K)t-Ks}-\e^{2( \|\mathcal R\|_\infty -K)s}}{2 \|\mathcal R\|_\infty-K}\d s\\
 &\ge P_tf^2-(P_tf)^2-\ff{\e^{Kt}-1}K  |\nn P_t f|^2.\end{align*}
 \item[$(4)$]  For any $f\in C_0^\infty(M)$ and $t\ge 0$,
\beg{align*} & P_tf^2-(P_tf)^2- \ff{1-\e^{-Kt}}K  P_t |\nn f|^2
& \ge - d P_t \|\Hess_{f}\|^2  \int_0^t \ff{\e^{2( \|\mathcal R\|_\infty-K)s}-\e^{-Ks}}{2 \|\mathcal R\|_\infty-K}\d s.\end{align*}\end{enumerate}
 \end{thm}
\beg{proof} By using
\beq\label{LOPP}\|\Hess_f\|_{HS}^2\le d \|\Hess_f\|^2,\end{equation}
the proof is completely similar to that of Theorem \ref{TP3}. For instance, below we only show the proof of (1) implying (2).

By $\Ric\le K$ and \eqref{QBC}, we have
$$\d |\nn P_{t-s}f|^2(X_s)\le  \big(K |\nn P_{t-s}f|^2 + \|\Hess_{P_{t-s}f}\|_{HS}^2 \big)(X_s)\d s+
2 \<\nn |\nn P_{t-s}f|^2(X_s), \Phi_s\d B_s\>$$   for $s\in [0,t].$
Combining this with \eqref{HS1} and \eqref{LOPP}, we arrive at
\beg{align*} &P_t |\nn f|^2 - \e^{Kt} |\nn P_t f|^2 \le d\int_0^t \e^{K(t-s)} P_s\|\Hess_{P_{t-s} f}\|^2\d s\\
&\le d\int_0^t \e^{K(t-s)} \e^{2( \|\mathcal R\|_\infty -K)(t-s)} P_t\|\Hess_{f}\|^2\d s\\
&= \ff{\e^{(2 \|\mathcal R\|_\infty-K)t}-1}{ 2 \|\mathcal R\|_\infty-K}    P_t\|\Hess_{f}\|^2 .\end{align*}
Then (1) implies (2). We therefore omit other proofs. \end{proof}

\section{Formula of   $\nn\Ric$}


 \beg{thm}\label{T1.2}    For any $x\in M$, $v_1,v_2\in T_xM$ and $f\in C_b^4(M)$ with $\nn f(x)=v_1, \Hess_f(x)=0$, there holds
 \beq\label{HD}  (\nn_{v_2}\Ric)(v_1,v_1) =2 \lim_{t\downarrow 0} \ff{(P_t\Hess_f-\Hess_{P_tf})(v_1,v_2) }{t}= \big(\DD\Hess_f -\Hess_{\DD f}\big)(v_1,v_2).  \end{equation}Consequently, $M$ is Ricci parallel if and only if
$\DD\Hess_f =\Hess_{\DD f}$ holds at any point $x\in M$ and $f\in C_0^\infty(M)$ with $\Hess_f(x)=0$.
  \end{thm}
When $f\in C_0^\infty(M)$, the second equation in \eqref{HD} follows from \eqref{DS5}. By a standard approximation argument, this equation holds for all $f\in C_b^2(M).$ So, it suffices to prove the first equation or
the formula \eqref{DRV} below for  $x\in M$ and $f\in C_0^\infty(M)$ with $\Hess_f(x)=0$. Here, we prove both of them by using analytic and probabilitstic arguments respectively, since each proof has its own interest.

\beg{proof}[Analytic Proof]   For any $x\in M$ and $f\in C_0^\infty(M)$ with $\Hess_f(x)=0,$ we intend to prove
\beq\label{DRV} (\nn_v \Ric)(\nn f,\nn f) =(\DD\Hess_f)(\nn f, v)- \Hess_{\DD f}(\nn f, v),\ \ v\in T_xM.\end{equation}
According to the Bochner-Weitzenb\"ock formula, we have
\beq\label{BW22}\Ric^\#(\nn f)= \DD\nn f- \nn \DD f,\end{equation}
where $\DD\nn f:= \Phi(\DD_{O(M)} \nn_f^{O(M)})(\Phi)$ is independent of $\Phi\in O(M).$ Consequently,
\beq\label{BW11} \Ric(\nn f,\nn f)= \ff 1 2\DD |\nn f|^2 -\<\nn \DD f,\nn f\>- \|\Hess_f\|_{HS}^2.\end{equation}
Since $\Hess_f(x)=0$, \eqref{BW11} and \eqref{BW22}  imply that  at point $x$,
\beg{align*}&(\nn\Ric)(\nn f,\nn f)= \nn\{\Ric(\nn f,\nn f)\} \\
&= \ff 1 2 \nn \DD |\nn f|^2 - \Hess_{\DD f}^\#(\nn f)-\Hess_f^\#(\nn\DD f)-2\|\Hess_f\|_{HS} \nn \|\Hess_f\|\\
&= \ff 1 2 \DD \{\nn |\nn f|^2\} -\ff 1 2 \Ric^\#(\nn |\nn f|^2) - \Hess_{\DD f}^\#(\nn f)\\
&= \DD \{\Hess_f^\#(\nn f)\}- \Hess_{\DD f}^\#(\nn f)\\
&= (\DD \Hess_f)^\#(\nn f) +\Hess_f^\#(\DD \nn f) + 2 {\rm tr}\big\{(\nn_\cdot \Hess_f^\#)(\Hess_f^\#(\cdot))\big\}- \Hess_{\DD f}^\#(\nn f)\\
&= (\DD \Hess_f)^\#(\nn f) - \Hess_{\DD f}^\#(\nn f).\end{align*}
Therefore, \eqref{DRV} holds. \end{proof}

\beg{proof}[Probabilistic Proof]  We first consider  bounded $\nn \Ric$ and $\mathcal R$, then extend to the general case by using Proposition \ref{P2}.

\paragraph{(a)} Assume that $\|\mathcal R\|_\infty +\|\nn\Ric\|_\infty<\infty.$ Let $x\in M$ and $v_1,v_2\in T_xM$. We take $f\in C_b^4(M)$ such that $\nn f(x)=v_1$ and $\Hess_f(x)=0$. Below we only  consider functions taking value at point $x$.  Since $\Hess_f(x)=0$, there exists a constant $c>0$ such that
\beq\label{OM}P_t\|\Hess_f\|_{HS}^2 (x)= \ff 1 2 \int_0^t P_s \DD \|\Hess_f\|_{HS}^2(x)\d s \le c t,\ \ t\ge 0.\end{equation} Then there exists a constant $c_1>0$ such that
\beq\label{M1} P_s \|\Hess_f\|\le \ss{P_s \|\Hess_f\|^2}\le c_1 \ss s,\ \ s\in [0,1].\end{equation}
Since $\nn f(x)=v_1$ and $\nn P_s f(x)$ is smooth in $s$,  this together with  \eqref{HS3} yields
\beq\label{M2}\E |\nn f(X_s)- \itparallel_s v_1|\le c_2   s,\ \ s\in [0,1]\end{equation} for some constant $c_2>0$.
Moreover, by \eqref{WT} there exists a constant $c_3>0$ such that
\beq\label{M3} | W_s(v_i)-\itparallel_s v_i|\le c_3 s, \ \ s\in [0,1], i=1,2.\end{equation}
Combining \eqref{M1}-\eqref{M3} with \eqref{AS3}, \eqref{W2}   and \eqref{W2'}, for small $t>0$ we arrive at
\beq\label{M5} \beg{split}  \E\<\nn f(X_t), W_t^{(2)} (v_1,v_2)\> = &\ff 1 2 \E  \int_0^t \big\<(\tt\nn \Ric^\#)(W_s(v_2))W_s(v_1), \nn f(X_s)\>\d s\\
& + \sum_{i=1}^d \E \int_0^t \Hess_{P_{t-s}f}\big(\Phi_s^i, {\mathcal R} (\Phi_s^i, W_s(v_2)) W_s(v_1)\big)\d s\\
= &\, {\rm o}(t) + \ff 12\<(\tt\nn \Ric^\#)(v_2)v_1,v_1\>t=  {\rm o}(t) - \ff 12(\nn_{v_2}\Ric)(v_1,v_1) t,\end{split}\end{equation}where the last step follows from \eqref{W2'}.

 On the other hand, by \eqref{HESS},  \eqref{M1} and \eqref{M3}, there exist a constant  $c_4 >0$ such that for small $t>0$,
 \beg{align*} & \E  \<\nn f(X_t), W_t^{(2)} (v_1,v_2)\>  =  \Hess_{P_tf}(v_1,v_2) -\E \Hess_f(W_t(v_1),W_t(v_2))  \\
 & =\Hess_{P_tf}(v_1,v_2) - P_t  \Hess_f (v_1,v_2)  + {\rm O}(t) P_t \|\Hess_f\| \\
 &=  \Hess_{P_tf}(v_1,v_2) - P_t  \Hess_f (v_1,v_2)  + {\rm o}(t).\end{align*}
 Combining this with \eqref{M5}   we derive the desired the first equation in \eqref{HD}.

 \paragraph{(b)} In general, let $M$ be a complete non-compact Riemannian manifold. For fixed $x\in M$ we take $D= B(x,4)$. Let $f\in C^\infty(\bar D)$ with $f|_{\pp D}=0, f|_{B(x,3)}=1$, $|\nn f|=1$ on $\pp D$ and $f>0$ in $D$.   Then $(D, f^{-2}g)$ is a complete Riemannian manifold. We use superscript $D$ to denote quantities on this manifold, for instance, $\mathcal R^D$ is the Riemannian tensor on $(D, f^{-2}g)$. Then both $\mathcal R^D$ and $\nn^D\Ric^D$ are bounded.   So, by (a), for $P_t^D$   the   heat semigroup on $D$,
 \beq\label{HSU1} (\nn_{v_2}\Ric)(v_1,v_1)=\lim_{t\downarrow 0} \ff{(P_t^D\Hess^D_f-\Hess^D_{P_t^Df})(v_1,v_2) }{t}.  \end{equation}
 Since $f=1$ in $B(x,3)$, we may construct the horizontal  Brownian motion $\Phi_t^D(y)$ on $D$ with $\Phi^D_0(y)=\Phi_0\in O_y(M)$  such that
 $$\Phi_t^D(y)= \Phi_t(y),\ \ t\le \tau_3(y),\ \ X_t^D(y)=X_t(y), t\le \tau_3(y),$$ where
 $$\tau_3(y)=\inf\{t\ge 0: X_t(y)\in\pp B(x,3)\}.$$
  Noting that
 $$P_t f(y)= \E [f(X_t(y))1_{\{t<\zeta\}}],\ \ P_t^Df(y)= \E[f(X_t^D(y))],\ \ f\in \B_b(M),$$ where $\zeta$ is the life time of $X_t$, we obtain
$$ |P_tf(y)-P_t^D(y)| \le \|f\|_\infty \P(\tau_3(y)\le t),\ \ t\ge 0, f\in \B_b(M).$$
Combining this with    \cite[Lemma 2.3]{ATW09}, we may find   constants $c_1, c_2>0$ such that
\beq\label{HSU2}|P_t f(y)- P_t^D f(y)|\le c_1 \|f\|_\infty \e^{-c_2/t},\ \ t\in (0,1], y\in B(x,2).\end{equation}
Consequently,
\beq\label{HSU3}|P_t\Hess_f(v_1,v_2)- P_t^D\Hess_f(v_1,v_2)|\le c_1 \|\Hess_f\|_\infty \e^{-c_2/t},\ \ t\in (0,1].\end{equation}
Moreover,   since $X_t=X_t^D\in B(x,2)$ before time $\tau_2$,   Proposition \ref{P2} and \eqref{HSU2} imply that at point $x$,
\beg{align*} &|\Hess_{P_tf}(v_1,v_2)- \Hess^D_{P_t^Df}(v_1,v_2)|\\
&\le \E\big[(|P_{t-t\land\tau_1}f(X_{t\land \tau_1})-P_{t-t\land\tau_1}^Df(X^D_{t\land\tau_1})|\cdot |M_t|\big]\\
&\qquad +\E\big[|P_{t-t\land\tau_2}f(X_{t\land \tau_2})-P_{t-t\land\tau_2}^Df(X^D_{t\land \tau_2})|) |N_{t}| \big]\\
&\le  c_1\e^{-c_2/t}\|f\|_\infty\ff{C(x)}{t},\ \ t\in (0,1]. \end{align*}
Combining this with \eqref{HSU1} and \eqref{HSU3}, we prove the first equation in \eqref{HD}.
 \end{proof}

\section{Identification of constant curvature}

\beg{thm}\label{TN} Let $k\in \R$. Then each of the following assertions is equivalent to $\Sect =k$: \beg{enumerate}
\item[$(1)$]  For any $ t\ge 0$ and $f\in C_0^\infty(M),$
\beq\label{SET} \Hess_{P_t f}= \e^{-d kt} P_t\Hess_f + \ff 1 d (1-\e^{-dkt}) (P_t\DD f)\g.\end{equation}
\item[$(2)$]  For any   $f\in C_0^\infty(M),$
\beq\label{SET2} \Hess_{\DD f}-\DD \Hess_f  = 2k (\DD f)\g-2dk \Hess_f.\end{equation}
\item[$(3)$]  For any  $x\in M, u\in T_xM$ and  $f\in C_0^\infty(M)$ with $\Hess_f(x)= u\otimes u\ ($i.e. $\Hess_f(v_1,v_2)= \<u,v_1\>\<u,v_2\>, v_1,v_2\in T_xM)$,
 \beq\label{CVS}\big(\Hess_{\DD f} - \DD\Hess_f\big)(v,v)= 2k \big(|u|^2|v|^2-\<u,v\>^2\big),\ \ v\in T_xM.\end{equation}
\item[$(4)$]   For any   $f\in C_0^\infty(M),$
\beq\label{CVS2}\ff 12 \DD \|\Hess_f\|_{HS}^2 -\<\Hess_{\DD f}, \Hess_f\>_{HS} -\|\nn \Hess_f\|_{HS}^2 = 2k \big(d\|\Hess_f\|_{HS}^2 -(\DD f)^2\big).\end{equation}\end{enumerate} \end{thm}

To prove this result, we need the following lemma where $\mathcal RT$ is defined in \eqref{RT}.

\beg{lem}\label{L*} If the sectional curvature $\Sect =k$ for some constant $k$, then for any symmetric $2$-tensor $T$,
$$\mathcal RT = k {\rm  tr}(T) \g -k T.$$
\end{lem}
\beg{proof}  Let $\tt T= k {\rm  tr}(T) \g -k T.$ Since both $\mathcal R T$ and $\tt T$ are symmetric, it suffices to prove
\beq\label{SE1} (\mathcal R T)(v,v)= \tt T(v,v),\ \ v\in T_x, |v|=1, x\in M.\end{equation}
Let $v\in T_xM$ with $|v|=1$. By the symmetry of $T$, there exists $\Phi=(\Phi^i)_{1\le i\le d} \in O_x(M) $ such that
$$T^\# (\Phi^i)= \ll_i \Phi^i,\ \ 1\le i\le d$$ holds for some constants $\ll_i, 1\le i\le d.$ Then $\Sect =k$ implies
\beg{align*} (\mathcal R T)(v,v) &= \sum_{i=1}^d \big\<\mathcal R(\Phi^i, v)v, T^\#(\Phi^i)\big\> =\sum_{i=1}^d \ll_i \big\<\mathcal R (\Phi^i, v)v,\Phi^i\big\>\\
&= k \sum_{i=1}^d \ll_i\big(1-\<\Phi^i, v\>^2\big) = k {\rm tr}(T)-k T(v,v)=\tt T(v,v).\end{align*} Therefore, \eqref{SE1} holds.
\end{proof}

\beg{proof}[Proof of Theorem \ref{TN}] Obviously, (3) follows from (2). Next, by \eqref{DS5}  and taking derivative of \eqref{SET} with respect to $t$ at $t=0$, we obtain \eqref{SET2}. So, (1) implies (2).
 Moreover, by chain rule we have
 \beq\label{ZJ} \ff 1 2 \DD\|\Hess_f\|_{HS}^2 =\<\DD\Hess_f,\Hess_f\>_{HS} +\|\nn \Hess_f\|_{HS}^2.\end{equation}
 Then (4) follows from (2) and the identity $\<(\dd f)\g, \Hess_f\>_{HS} =(\DD f)^2.$ To complete the proof,
 below we prove $``\Sect=k \Rightarrow (1)$", $``(3) \Rightarrow\Sect =k$ " and $``(4) \Rightarrow\Sect =k$ "   respectively.

{\bf (a)} $\Sect=k \Rightarrow (1)$.  Let $\Sect =k$. Then $\Ric= (d-1)k$.  By \eqref{WT}, \eqref{HS5''}, we have $W_s(v)= \e^{-\ff K 2 s} \parallel_s v,\ s\ge 0, v\in TM$, and for any $x\in M, v_1,v_2\in T_xM$,
$$ \Hess_{P_tf} (v_1,v_2)  = \e^{-Kt} P_t \Hess_f (v_1,v_2) +\int_0^t \e^{-Ks} P_s (\mathcal R \Hess_{P_{t-s}f})(v_1,v_2)\d s.$$
 Noting that $\DD P_{t-s}f = P_{t-s}\DD f$ and $\<\parallel_s v_1, \parallel_s v_2\>=\<v_1,v_2\>$, this together with    Lemma \ref{L*} gives
\beg{align*} &\Hess_{P_tf} (v_1,v_2) - \e^{-Kt} P_t \Hess_f(v_1,v_2)\\
&= \int_0^t \e^{-Ks}\Big(P_s\big\{k\<v_1,v_2\> P_{t-s}\DD f\big\} -k P_s(\Hess_{P_{t-s}f})(v_1,v_2)\Big)\d s,\ \ t\ge 0.\end{align*}  Therefore,
\beg{align*} \Hess_{P_tf} (v_1,v_2) &= \e^{-(K+k)t} P_t \Hess_f(v_1,v_2) + k\<v_1,v_2\>\int_0^t \e^{-(K+k)s} P_t \DD f \d s \\
&= \e^{-dkt}P_t \Hess_f(v_1,v_2) +\ff{1-\e^{-dkt}}d (P_t\DD f) \<v_1,v_2\>.\end{align*} So, \eqref{SET} holds.

{\bf (b)} $(3) \Rightarrow\Sect =k$. By taking $u=0$, (3) implies that for any $f\in C_0^\infty(M)$ with $\Hess_f(x)=0$,
$$(\DD\Hess_f-\Hess_{\DD f})(v,v)=0,\ \ v\in T_xM.$$ By the symmetry of $\DD\Hess_f-\Hess_{\DD f}$, this is equivalent to
$$(\DD\Hess_f-\Hess_{\DD f})(v_1,v_2)=0,\ \ v_1,v_2\in T_xM.$$ So,
 for any $v_1,v_2\in T_xM$,  by taking
 $f\in C_0^\infty(M)$ such that $\nn f(x)= v_1$ and $\Hess_f(x)=0$, we deduce from Theorem \ref{T1.2} that
$$(\nn_{v_2}\Ric)(v_1,v_1)=  \big(\DD\Hess_f - \Hess_{\DD f}\big)(v_1,v_2)= 0.$$ Thus, $M$ is Ricci parallel.  By Theorem \ref{T3.1}, \eqref{HS5''} holds.
Due to \eqref{WT} and It\^o's formula, by taking derivative of \eqref{HS5''} with respect to $t$ at $t=0$, we obtain
\beq\label{DH5} \beg{split}&\ff 1 2 \Hess_{\DD f} (v_1,v_2)\\
&=   \ff 1 2 (\DD \Hess_f)(v_1,v_2) -  \Ric(v_1, \Hess_f^\#(v_2))
  + {\rm tr} \<\mathcal R(\cdot, v_2)v_1, \Hess_f^\#(\cdot)\>,\ \ f\in C_0^\infty(M). \end{split}\end{equation}
Now, letting $u,v\in T_xM$ with $|u|=|v|=1$ and $\<u,v\>=0$, and combining \eqref{DH5}  with \eqref{CVS} for $v_1=v_2=v$,   and  $f\in C_0^\infty(M)$ with $\Hess_f(x)= u\otimes u$, we arrive at
\beg{align*}  k&=k(|u|^2|v|^2-\<u,v\>^2) = \ff 1 2 (\Hess_{\DD f}- \DD\Hess_f)(v,v)   \\
&= - \Ric(v,\Hess_f^\#(v)) +  {\rm tr} \<\mathcal R(\cdot, v)v, \Hess_f^\#(\cdot)\>= \Sect(u,v).\end{align*}Therefore, $\Sect =k.$

{\bf (c)} $(4) \Rightarrow\Sect =k$. By \eqref{ZJ}, \eqref{CVS2} is equivalent to
\beq\label{ZJ2} \ff 1 2 \big\<\DD\Hess_f- \Hess_{\DD f}, \Hess_f\big\>_{HS} =k \big(d \|\Hess_f\|_{HS}^2-(\DD f)^2\big).\end{equation}

We first prove that \eqref{ZJ2} implies $\nn \Ric=0$. Let $f\in C_0^\infty(M)$ with $\Hess_f(x)=0$. For any $u\in T_xM$, let  $h\in C_0^\infty (M)$ with $\Hess_h(x)=u\otimes u$.
Applying \eqref{ZJ2} to $f_s:= f+ sh, s\ge 0$, we obtain at point $x$ that
\beg{align*} &\ff s 2 (\DD\Hess_f- \Hess_{\DD f})(u,u) +\ff{s^2} 2 (\DD\Hess_h-\Hess_{\DD h})(u,u)\\
&= k \big(ds^2\|\Hess_h\|_{HS}^2 -s^2 (\DD h)^2\big),\ \ s>0.\end{align*}
Multiplying by $s^{-1}$ and letting $s\to 0$, we arrive at $(\DD\Hess_f-\Hess_{\DD f})(u,u)=0.$ As shown above, this implies $\nn \Ric=0$.

Next, we prove $\Sect=k$. Since $\nn\Ric=0$, \eqref{DH5} holds. 
 For $x\in T_xM$ and $u,v\in T_xM$ with $|u|=|v|=1$ and $\<u,v\>=0$, take $f\in C_0^\infty(M)$ such that
$$\Hess_f(x)= u\otimes v+ v\otimes u.$$ Then at point $x$,
$$\Hess_f^\#(\cdot)= \<u,\cdot\>v+ \<v,\cdot\>u.$$ So, by \eqref{DH5} and \eqref{ZJ2} we obtain
\beq\label{POL} \beg{split}  &2kd = k\big(d \|\Hess_f\|_{HS}^2 -(\DD f)^2\big)(x)\\
& =\ff 1 2 \< \DD\Hess_f- \Hess_{\DD f}, \Hess_f\>_{HS}(x) = (\DD\Hess_f -\Hess_{\DD f})(u,v) \\
&= 2 \Ric (u, \Hess_f^\#(v)) -2{\rm tr} \<\mathcal R(\cdot, v)u, \Hess_f^\#(\cdot)\> = 2 \Ric(u,u) + 2 \Sect(u,v).\end{split}\end{equation}
Letting $\{v_1\}_{1\le i\le d-1}$ be orthonormal and orthogonal to $u$, replacing $v$ by $v_i$ and sum over $i$ leads to
$$2kd(d-1) = 2(d-1)\Ric(u,u)+2\Ric(u,u)= 2d\Ric(u,u).$$
Thus, $\Ric(u,u)= (d-1)k$, and \eqref{POL} implies $\Sect(u,v)= k.$

\end{proof}

\section{Identifications of Einstein   manifolds}

\beg{thm}\label{T1.1} For any constant $K\in \R$, the following statements are equivalent each other:
\beg{enumerate} \item[$(1)$] $M$ is an Einstein manifold with $\Ric= K$.
\item[$(2)$]  $ \|\mathcal R\|_\infty <\infty$, and for any $x\in M, t\ge 0$, $f\in C_0(M)$ and $v_1,v_2\in T_xM$,
 \beq\label{HS50} \Hess_{P_tf}(v_1,v_2) - \e^{-Kt} P_t\Hess_f(v_1, v_2)
 = \int_0^t\e^{-Ks}P_s\big(\mathcal R \Hess_{P_{t-s}f}\big)  (v_1, v_2) \d s. \end{equation}
 \item[$(3)$] $\|\mathcal R\|_\infty<\infty,$ and
\beg{align*} &\ff{\e^{Kt}- \e^{2(K- \|\mathcal R\|_\infty)t}}{ 2 \|\mathcal R\|_\infty-K}  \|\Hess_{P_tf}\|_{HS}^2  \le P_t|\nn f|^2 - \e^{Kt} |\nn P_t f|^2\\
& \le \ff{\e^{(2 \|\mathcal R\|_\infty-K)t}-1}{2(d -1)\|\mathcal R\|_\infty-K} P_t \|\Hess_f\|^2_{HS},\ \  f\in C_b^2(M), t\ge 0.\end{align*}
  \item[$(4)$] There exists   $h:  [0,\infty)\times M\to   [0,\infty)$ with $\lim_{t\to 0} h(t,\cdot)=0$ such that
   $$ \big|P_t|\nn f|^2 - \e^{Kt} |\nn P_t f|^2 \big|\le h(t,\cdot) \big(\|\Hess_{P_tf}\|_{HS}^2+P_t\|\Hess_f\|_{HS}^2\big),\ \ t\ge 0, f\in C_0^\infty(M).$$
   \item[$(5)$] $\|\mathcal R\|_\infty<\infty,$ and
\beg{align*} &\|\Hess_{P_tf}\|_{HS}^2\int_0^t \ff{\e^{Ks} -\e^{2(K- \|\mathcal R\|_\infty)s}}{2 \|\mathcal R\|_\infty-K}\d s
 \le P_tf^2-(P_tf)^2-\ff{\e^{Kt}-1}K  |\nn P_t f|^2\\
& \le \big(P_t \|\Hess_f\|^2_{HS}\big) \int_0^t \ff{\e^{(2 \|\mathcal R\|_\infty -K)t - Ks}- \e^{2( \|\mathcal R\|_\infty-K)s}}{2  \|\mathcal R\|_\infty-K}\d s,\ \ f\in C_b^2(M), t\ge 0.\end{align*}
\item[$(6)$] $\|\mathcal R\|_\infty<\infty,$ and
\beg{align*} &-\big(P_t\|\Hess_f\|^2_{HS} \big) \int_0^t \ff{\e^{2( \|\mathcal R\|_\infty-K)s}-\e^{-Ks})}{2  \|\mathcal R\|_\infty-K}\d s
 \le P_tf^2-(P_tf)^2-\ff{1-\e^{-Kt}}K  P_t |\nn f|^2\\
& \le - \|\Hess_{P_tf}\|^2_{HS}\,\e^{2(K- \|\mathcal R\|_\infty)t} \int_0^t \ff{ \e^{2( \|\mathcal R\|_\infty-K)s}-\e^{-Ks}}{2 \|\mathcal R\|_\infty-K}\d s, \ \ f\in C_b^2(M), t\ge 0.\end{align*}
\item[$(7)$] There exists   $\tt h:  [0,\infty)\times M\to   [0,\infty)$ with $\lim_{t\to 0}t^{-1} \tt h(t,\cdot)=0$ such that
\beg{align*} &\min\bigg\{\Big|P_tf^2-(P_tf)^2-\ff{\e^{Kt}-1}K  |\nn P_t f|^2\Big|,\ \Big|P_tf^2-(P_tf)^2-\ff{1-\e^{-Kt}}K  P_t |\nn  f|^2\Big|\bigg\}\\
 &\le     \tt h(t,\cdot) \big(\|\Hess_{P_tf}\|_{HS}^2+ P_t  \|\Hess_f\|^2_{HS}\big),\ \ t\ge 0, f\in C_0^\infty(M).\end{align*}
 \item[$(8)$] For any $f\in C_0^\infty(M)$,
 $$\ff 1 2 \big\{\Hess_{\DD f}- \DD\Hess_f\big\}= (\mathcal R\Hess_f) - K\Hess_f.$$
\end{enumerate}\end{thm}

\beg{proof} Obviously,     (3)  implies (4), each of (5) and (6) implies (7). According to Theorem \ref{T3.1} for $\Ric=K$, (1) implies (2). Moreover, by taking derivative of \eqref{HS50} with respect to $t$ at $t=0$, we obtain (8).
So, it suffices to  prove that   (1) implies (3);  (3) implies (5) and  (6); and  each of     (4), (7) and (8) implies (1).

\paragraph{(a)} $(1)\Rightarrow (3)$.
By the Bochner-Weitzenb\"ock formula and using $\Ric =K$,  we obtain
\beg{align*} &\ff 1 2 \DD |\nn P_{t-s}f|^2(X_s)-\<\nn P_{t-s}f, \nn \DD P_{t-s}f\>(X_s)\\&= \Ric(\nn P_{t-s}f, \nn P_{t-s}f)(X_s) +\|\Hess_{P_{t-s}f}\|_{HS}^2(X_s)\\
&=K |\nn P_{t-s}f|^2(X_s) + \|\Hess_{P_{t-s}f}\|_{HS}^2(X_s).\end{align*} Then \eqref{QBC} implies
$$\d |\nn P_{t-s}f|^2(X_s) = \big(K |\nn P_{t-s}f|^2 + \|\Hess_{P_{t-s}f}\|_{HS}^2 \big)(X_s)\d s+
2 \<\nn |\nn P_{t-s}f|^2(X_s), \Phi_s\d B_s\>$$   for $s\in [0,t].$
Thus,
\beq\label{A1} P_t |\nn f|^2 - \e^{Kt} P_t|\nn f|^2 =\int_0^t \e^{K(t-s)} P_s\|\Hess_{P_{t-s} f}\|_{HS}^2\d s.\end{equation}
Since by \eqref{HS2}
$$\|\Hess_{P_{t-s}f}\|_{HS}^2 \le \e^{2( \|\mathcal R\|_\infty -K)(t-s)} P_{t-s}\|\Hess_f\|_{HS}^2,$$
it follows from \eqref{A1} that
\beg{align*} P_t |\nn f|^2 - \e^{Kt} P_t|\nn f|^2&\le \int_0^t \e^{K(t-s) + 2( \|\mathcal R\|_\infty -K)(t-s)} P_t\|\Hess_f\|_{HS}^2\d s\\
&= \ff{\e^{(2 \|\mathcal R\|_\infty-K)t}-1}{2  \|\mathcal R\|_\infty-K} P_t \|\Hess_f\|^2.\end{align*}
So, the second inequality in (3) holds. Similarly, \eqref{HS2} implies
$$P_s\|\Hess_{P_{t-s}f}\|_{HS}^2 \ge \e^{-2( \|\mathcal R\|_\infty -K)s}  \|\Hess_{P_sP_{t-s}f}\|_{HS}^2 =\e^{-2( \|\mathcal R\|_\infty -K)s}  \|\Hess_{P_tf}\|_{HS}^2,$$
the first inequality in (3) also follows from \eqref{A1}.

\paragraph{(b)} $(3)\Rightarrow (5)$ and $(6)$.
By (3) and \eqref{HS2} we have
\beg{align*} &\ff{\e^{Ks}-\e^{2(K- \|\mathcal R\|_\infty)s}}{2 \|\mathcal R\|_\infty-K}\|\Hess_{P_tf}\|_{HS}^2\le P_s|\nn P_{t-s}f|^2 -\e^{Ks}|\nn P_tf|^2\\
&\le \ff{\e^{(2  \|\mathcal R\|_\infty-K)s}-1}{2 \|\mathcal R\|_\infty-K}P_s\|\Hess_{P_{t-s}f}\|_{HS}^2\\
&\le \ff{\e^{((2  \|\mathcal R\|_\infty-K)t-K(t-s)}-\e^{2( \|\mathcal R\|_\infty-K)(t-s)}}{2\|\mathcal R\|_\infty-K}P_t\|\Hess_f\|_{HS}^2.\end{align*}
This together with \eqref{PTT} ensures (5).

Similarly, \eqref{HS2} and (3) imply
\beg{align*} & \ff{\e^{-Ks}(\e^{(2 \|\mathcal R\|_\infty-K)s}-1)}{2 \|\mathcal R\|_\infty-K}P_t\|\Hess_f\|_{HS}^2\ge \e^{-Ks} P_t|\nn f|^2- P_{t-s}|\nn P_s f|^2\\
&\ge \ff{1-\e^{(K-2  \|\mathcal R\|_\infty)s}}{2 \|\mathcal R\|_\infty-K} P_{t-s}\|\Hess_{P_sf}\|_{HS}^2\\
&\ge \e^{2(K- \|\mathcal R\|_\infty)t} \|\Hess_{P_tf}\|_{HS}^2 \cdot \ff{\e^{2( \|\mathcal R\|_\infty-K)s}-\e^{-Ks}}{2 \|\mathcal R\|_\infty-K},
\end{align*}  which together with \eqref{PTT} gives (6).

\paragraph{(c)} $(4)\Rightarrow (1)$. For  $v\in T_xM$ with $|v|=1$, take $f\in C_0^\infty(M)$ such that
$$\nn f(x)=v,\ \ \Hess_f(x)=0.$$ Then \eqref{BW} holds. Moreover, since $\Hess_{P_tf}(x)$ is smooth in $t$, $\Hess_f(x)=0$ implies
$$\|\Hess_{P_t f}(x)\|\le c(x) t,\ \ t\in [0,1]$$ for some constant $c(x)>0$.
Combining this with   \eqref{BW}, \eqref{OM} and  (4), we obtain
\beg{align*} &0= \lim_{t\downarrow 0} \ff{ P_t|\nn f|^2(x)-\e^{Kt}|\nn P_t f|^2(x)}{t} \\
&= \lim_{t\downarrow 0} \bigg( \Ric(v,v) +\ff{1-\e^{Kt}}{t} P_t|\nn f|^2(x)\bigg)=    \Ric(v,v)-K.\end{align*}
Therefore, $\Ric(v,v)=K|v|^2$ holds for all $v\in TM$. By the symmetry of $\Ric$ and $g$, this is equivalent to   $\Ric=K$.

\paragraph{(d)} $(7)\Rightarrow (1)$. Let $v$ and $f$ be in {\bf (c)}. In the spirit of  \eqref{BW} and using \eqref{PTT},  we have
\beg{align*} &2\lim_{t\downarrow 0} \ff {P_t f^2-(P_tf)^2 -\ff{\e^{Kt}-1}K |\nn P_tf|^2} {t^2} \\
&= \lim_{t\downarrow 0} \ff {2}{t^2} \int_0^t \big\{P_s|\nn P_{t-s}f|^2- \e^{Ks} |\nn P_t f|^2\big\}(x)\d s= \Ric(v,v)-K,\end{align*}
and
\beg{align*} &2\lim_{t\downarrow 0} \ff {P_t f^2-(P_tf)^2 -\ff{1-\e^{-Kt}}K P_t |\nn f|^2} {t^2} \\
&= \lim_{t\downarrow 0} \ff {2}{t^2} \int_0^t \big\{P_s|\nn P_{t-s}f|^2- \e^{-Ks} P_t |\nn  f|^2\big\}(x)\d s= K-\Ric(v,v).\end{align*}
Thus, multiplying  the inequality  in (7) by $t^{-2}$ and letting $t\to 0$, we prove $ \Ric(v,v)-K =0.$ That is, (1) holds.

\paragraph{(e)} $(8)\Rightarrow (1)$.  For any $v_1,v_2\in T_xM$, take $f\in C_0^\infty(M)$ such that $\nn f(x)= v_1, \Hess_f(x)=0$. According to Theorem \ref{T1.2}, (8) implies
$$(\nn_{v_2}\Ric)(v_1,v_1) = 0.$$ So, $M$ is Ricci parallel, and as shown in the proof of Theorem \ref{T3.1} that \eqref{DH5} holds. Taking $v_1=v_2=v$ for $v\in T_xM$ with $|v|=1$, and letting $f\in C_0^\infty(M)$ such that $\Hess_f(x)= v\otimes v,$  \eqref{DH5} implies
$$\ff 1 2 \big\{\Hess_{\DD f}-\DD\Hess_f\big\}(v,v)= -\Ric(v,v)+{\rm tr}\<\mathcal R(\cdot, v)v, \Hess_f^\#(\cdot)\> =-\Ric(v,v).$$ Combining this with (8) we obtain
$$-\Ric(v,v)= -K\Hess_f(v,v)=-K.$$ So, (1) holds.

\end{proof}

\section{Identifications of Ricci Parallel manifolds}

\beg{thm}\label{T1.3} The following assertions are equivalent each other: \beg{enumerate} \item[$(1)$]
$M$ is a Ricci parallel manifold.
 \item[$(2)$]   $ \|\mathcal R\|_\infty <\infty$, and $\eqref{HS5''}$ holds for any $x\in M, t\ge 0$, $f\in C_0(M), v_1,v_2\in T_xM$.

 \item[$(3)$]  $ \|\mathcal R\|_\infty <\infty$, and for any constant $K\in\R$ with $\Ric\ge K$,  $t\ge 0, f\in C_b^2(M)$,
 \beq\label{LEO2}\|\Hess_{P_tf} -   P_t\Hess_f\|\le  \bigg(\ff{\|\Ric\|_\infty(1-\e^{-Kt})}K+  \e^{( \|\mathcal R\|_\infty -K)t}-\e^{-Kt}\bigg)P_t\|\Hess_f\|.\end{equation}
 \item[$(4)$]   There exists a function
 $h: [0,\infty)\times M \to  [0,\infty) $ with $\lim_{t\downarrow 0}t^{-\ff 1 2}h(t,\cdot)=0$ such that
\beq\label{LEO} \|\Hess_{P_tf} -   P_t\Hess_f\|\le h(t,\cdot)\big( P_t\|\Hess_f\|+\|\Hess_{P_tf}\|\big),\ \ t\ge 0, f\in C_0^\infty(M).\end{equation}
 \item[$(5)$] For any $f\in C_0^\infty(M)$ and $x\in M$,
  $$\big(\Hess_{\DD f}- \DD \Hess_f\big)(v_1,v_2)= 2\big(\mathcal R\Hess_f\big)(v_1,v_2) - 2\Ric(v_1,\Hess_f^\#(v_2)),\ \ v_1,v_2\in T_xM.$$
\item[$(6)$] For any $x\in M$ and $f\in C_0^\infty(M)$ with $\Hess_f(x)=0$,
$$(\DD\Hess_f)(v_1,v_2)=\Hess_{\DD f}(v_1,v_2),\ \  v_1,v_2\in T_xM.$$
 \end{enumerate}
 \end{thm}

\beg{proof} The equivalence of (1) and (6) follows from Theorem \ref{T1.2}, (1) implying (2) is included in Theorem \ref{T3.1}, (5) follows from (2) by taking derivative of $\eqref{HS5''}$ with respect to $t$ at $t=0$,
and it is obvious that (3) implies (4) while   (6) follows from  (5).     So,  it remains to prove that $(1)$  implies $(3)$,    and  $(4)$ implies $(1)$.

\paragraph{(a)} $(1)$ implies $(3)$. Let $ M$ be Ricci parallel with $\Ric\ge K$.  By \eqref{WT}  we have
\beq\label{KL} |W_s(v)|\le \e^{-\ff K 2 s} |v|,\ \ v\in T_xM,\end{equation} and
\beg{align*}
&\d |W_t(v)-\itparallel_t v|^2= \d |\Phi_t^{-1}W_t(v)-\Phi_0^{-1} v|^2\\
&=\<W_t(v)-\itparallel_tv, \Ric^\#(W_t(v))\>\d t\le |W_t(v)-\itparallel_t v| \cdot\|\Ric\|_\infty\e^{-\ff K 2 t}|v|.\end{align*}
So,
$$|W_t(v)- v|\le \ff {\|\Ric\|_\infty}2 \int_0^t \e^{-\ff K 2 s}\d s= \ff{\|\Ric\|_\infty (1-\e^{-\ff K 2 t})}K,\ \ |v|\le 1.$$
Thus, for $|v_1|, |v_2|=1$,
\beq\label{TJ1}\beg{split}& \big|\Hess_f(W_t(v_1), W_t(v_2))- \Hess_f(\itparallel_t v_1, \itparallel_t v_2)\big| \\
 &\le \|\Hess_f\|\big(|W_t(v_1)|\cdot |W_t(v_2)-\itparallel_t v_2)| + |W_t(v_1)-\itparallel_t v_1|\big)\\
 &\le \|\Hess_f\|\big(\e^{-\ff K 2 t} +1\big)\big(1-\e^{-\ff K 2 t}\big) \ff{\|\Ric\|_\infty} K= \|\Hess_f\|\ff{\|\Ric\|_\infty(1-\e^{-Kt})}K.\end{split}\end{equation}
  Combining \eqref{HS5''}  with  \eqref{ZJJ} for $\zeta=\infty$,    \eqref{HS1}, \eqref{KL} and \eqref{TJ1}, we obtain
\beg{align*} &\|\Hess_{P_tf} - P_t\Hess_f\|\\
&\le \sup_{|v_1|,|v_2|\le 1}\big|\Hess_{P_tf} (v_1,v_2) - \E  \Hess_f(W_t(v_1),W_t(v_2))\big| + \ff{\|\Ric\|_\infty(1-\e^{-Kt})}K P_t\|\Hess_f\|\\
 &\le   \ff{\|\Ric\|_\infty(1-\e^{-Kt})}K  P_t\|\Hess_f\| +  \|\mathcal R\|_\infty \int_0^t \e^{-Ks} \,\E\|\Hess_{P_{t-s}f}\|(X_s)\d s\\
&\le  \ff{\|\Ric\|_\infty(1-\e^{-Kt})}K  P_t\|\Hess_f\| +  \|\mathcal R\|_\infty \int_0^t \e^{-Ks}P_s\|\Hess_{P_{t-s}f}\| \d s\\
&\le   P_t\|\Hess_f\|\bigg(\ff{\|\Ric\|_\infty(1-\e^{-Kt})}K  +  \|\mathcal R\|_\infty \int_0^t \e^{-Ks+( \|\mathcal R\|_\infty -K)(t-s)}\d s\bigg) \\
&= \Big(\ff{\|\Ric\|_\infty(1-\e^{-Kt})}K+  \e^{( \|\mathcal R\|_\infty-K) t}-\e^{-Kt}\Big) P_t\|\Hess_f\|.\end{align*} So,   \eqref{LEO2} holds.

\paragraph{(b)} $(4) \Rightarrow  (1)$.  For any $x\in M$ and $v_1,v_2\in T_xM$, let $f\in C_0^\infty(M)$ such that $\nn f(x)=v_1$ and $\Hess_f(x)=0$. Since $\Hess_{f}(x)=0$ and $\Hess_{P_tf}(x)$ is   smooth in $t\ge 0$,
 $$\|\Hess_{P_tf}\|(x)\le c t,\ \ t\in [0,1]$$ holds for some constant $c>0$. Combining this with  Theorem \ref{T1.2}, \eqref{LEO} and \eqref{OM}, we obtain
\beg{align*}& \big|(\nn_{v_2}\Ric)(v_1,v_1)\big| = \lim_{t\downarrow 0} \ff {|P_t\Hess_f-\Hess_{P_tf}|(v_1,v_2)}t \\
&\le \lim_{t\downarrow 0} \ff{h(t,x)} t \big(P_t\Hess_f\|(x)+\|\Hess_{P_tf}\|(x) \big) \le \lim_{t\downarrow 0} \ff{h(t,x)(c_1\ss t+ ct)}{ t}=0.\end{align*}
By the symmetry of $(\nn_{v_2}\Ric)$, this implies $\nn\Ric=0$. Thus, (4) implies (1).
\end{proof}


\end{document}